\documentclass[review]{elsarticle}

\usepackage{lineno,hyperref}
\modulolinenumbers[5]

\usepackage{geometry}
\usepackage{amssymb,amsmath,amssymb,mathrsfs}

\newtheorem{remark}{Remark}[section]
\newproof{proof}{Proof}
\numberwithin{equation}{section}

\usepackage{color}
\definecolor{red}{rgb}{1,0,0} 
\definecolor{blue}{rgb}{0,0,1}
\definecolor{darkblue}{rgb}{0,0,0.8}

\usepackage{graphicx}

\begin{document}
	
	\begin{frontmatter}
		
		\title{Asymptotic-preserving schemes for kinetic-fluid modeling of mixture flows with distinct particle sizes}	
		
		\author[mymainaddress,mysecondaryaddress,mythirdaddress]{Shi Jin}
		\ead{shijin-m@sjtu.edu.cn}
		
		\address[mymainaddress]{School of Mathematical Sciences, Shanghai Jiao Tong University, Shanghai200240, China}
		\address[mysecondaryaddress]{Institute of Natural Sciences, Shanghai Jiao Tong University, Shanghai200240, China}
		\address[mythirdaddress]{Ministry of Education, Key Laboratory in Scientific and Engineering Computing, Shanghai Jiao Tong University, Shanghai200240, China}
		
		\author[mymainaddress,mysecondaryaddress]{Yiwen Lin\corref{mycorrespondingauthor}}
		\cortext[mycorrespondingauthor]{Corresponding author}
		\ead{linyiwen@sjtu.edu.cn}

		\begin{abstract}
			We consider coupled models for particulate flows, where the disperse phase is made of particles with distinct sizes. We are thus led to  a system coupling the incompressible Navier–Stokes equations to the multi-component Vlasov–Fokker–Planck equations. We design an {\it asymptotic-preserving} numerical scheme to approximate the system. 
			The scheme is based on  suitable  implicit treatment of the stiff drag force term as well as the Fokker–Planck operator, and can be formally shown to capture the hydrodynamic limit with time step and mesh size independent of the Stokes number.
			Numerical examples illustrate the accuracy and asymptotic behavior of the scheme, with several interesting applications.
			
		\end{abstract}
		
		
	\end{frontmatter}
	
	
	\section{Introduction}
	
	This paper concerns the kinetic-fluid models for a mixture of flows in which the particles represent the disperse phase evolving in a dense fluid. 
	Applications of such kinetic-fluid models include the dispersion of smoke or dust \cite{Friedlander1977}, biomedical modeling of spray \cite{Baranger2005}, coupled models in combustion theory \cite{Williams1985}, etc. 
	Specifically, we focus on the models that describe a large number of particles, {\it with distinct but fixed sizes}, interacting with a fluid.  Here the dense fluid phase is modeled by the Euler or Navier-Stokes equations and particles dispersed in the fluid are modeled by Fokker-Planck type kinetic equations.
	Such multi-size particle systems have a wide range of applications in engineering, especially for the complex meteorological simulation of large aircraft icing process. 
	For the distribution of droplets in the air, the influence of large droplets contained in the droplet distribution cannot be ignored \cite{Cober2006,Hauf2006,Potapczuk2019}.
	However, it is very difficult to simulate multi-size particles by experimental means. 
	See \cite{Prosperetti2007, Goudon2013fluid} for more details on the modeling of such multi-phase flows. 
	
	In this paper, we study the model where the evolution of the particle distribution function is driven by a combination of particle transport, a drag force exerted by the surrounding fluid on the particle obeying the Stokes Law, an external force field (such as gravity, electrostatic force, centrifugal force, etc.) and Brownian motion of particles. 
	For the two physically important regimes first investigated by Goudon et al. \cite{GoudonJabin2004a, GoudonJabin2004b}, we focus on the fine particle regime given in \cite{GoudonJabin2004b}.
	The fluid phase is incompressible and viscous,
	all phases are isothermal, 
	with interactions including coagulation and fragmentation that occur between particles, while  the change of particle sizes is ignored.  
	For the sake of simplicity, we suppose that the fluid density is constant and homogeneous.
	In this model, particles and fluid systems are coupled through nonlinear terms. Such a coupling and nonlinearities pose new difficulties in mathematical analysis and numerical computations when compared with uncoupled problems.  
	Furthermore, from a numerical point of view, the kinetic framework leads to high computational costs in both size and time, posing further computational challenges.
	
	The study of existence, uniqueness and regularity problems depends on the nature of the coupling and the complexity of the equations used to describe the fluid. 
	For the two-phase flow model system, it is worth mentioning related works like existence of strong solutions locally in time without velocity-diffusion \cite{Baranger2006}, existence of weak solutions for the Vlasov–Stokes system \cite{Hamdache1998} and for the incompressible Vlasov–Navier–Stokes system on a periodic domain \cite{BoudinDesvillettes2009} or a bounded domain \cite{Yu2013}, global-in-time existence of classical solutions close to the equilibrium for the incompressible Navier–Stokes–Vlasov–Fokker–Planck system \cite{GoudonHe2010}, analysis of compressible models \cite{Mellet2007}, several studies of coupling with the Euler system without viscosity \cite{Carrilo2006, CarrilloDuan2011} and systems with energy exchanges \cite{Boudin2009fluid}.
	Analysis of the asymptotics in the two-phase flow system is due to \cite{ GoudonJabin2004a, GoudonJabin2004b} by means of relative entropy methods, see also \cite{Mellet2008}. 
	For the multi-phase flow model system \eqref{ModelEquation0}, existence has been discussed in \cite{Goudon2013fluid} and regularity properties of the solutions close to the equilibrium as well as its long time behavior have been investigated recently by the authors in \cite{JinLin2022}.
	
	Numerical methods for such particulate flows have been developed in recent years, including particle-in-cell method \cite{Andrews1996}, Eulerian–Lagrangian method \cite{ Patankar2001IJMFa, Patankar2001IJMFb}, level set approach \cite{ Liu2011} and so on.
	One of the difficulties in numerically solving such multi-component Vlasov-Fokker-Planck-Navier-Stokes systems comes from the varying  Stokes number $\varepsilon$, which describes the ratio of the Stokes settling time over a certain time unit of observation.
	Due to the multiscale nature of the problem, it is often desired to design numerical schemes that possess the asymptotic-preserving (AP) properties \cite{Carrillo2008, Goudon2012, Goudon2013asymptotic}.
	The AP schemes (coined in \cite{Jin1999}) refer to those that, when letting the Stoeks number $\varepsilon$ go to zero and holding the mesh size and time step fixed, the numerical schemes for coupled kinetic-fluid models automatically become good numerical schemes for the hydrodynamic limiting equations, with numerical stability independent of $\varepsilon$. 
	We refer to \cite{Jin2010, Jin2022} and references therein for  reviews on  AP schemes and their applications.
	Most of these references do not address the effect of different particle sizes, namely all particles in the model are assumed to be in the same size. 
	This paper considers kinetic-fluid models for a mixture of the flows for particles with {\it {distinct}} sizes.

	The goal of this work is to design a numerical scheme to simulate the behavior of the fluid-particles systems with disparate particle sizes, capable of handling different regimes, from $\varepsilon=O(1)$ (the kinetic regime) to $\varepsilon \ll 1$ (the hydrodynamic regime). We will follow the discretization introduced in \cite{JinYan2011}. Specifically, 
	the AP scheme for this multi-phase model uses a combination of the projection method \cite{Chorin1967, Chorin1969} for the Navier-Stokes equations and an implicit treatment for stiff Fokker-Planck operators. 
	In addition to the challenge of nonlinear coupling between particles and the fluid,  new difficulty arises here due to the multi-phase properties where the disperse phase is made of particles with \textit{distinct} sizes. Roughly speaking,
		one needs to prove the velocities of different species equilibrate which is a property one does not encounter for two-phase kinetic-fluid systems with identical particle size, and thus one has to investigate  different convergence rates to equilibria for particles with different sizes and to justify AP properties of the scheme under the multi-phase kinetic-fluid system with disparate particle sizes. 
		For the coupling and nonlinearities, the construction of the scheme relies on evaluating implicitly the stiff terms of the multi-phase system. 
		It will require a carefully-designed time splitting which allows to compute implicitly the stiff drag force term efficiently and an inversion of the Fokker-Planck operator. 
	In the hydrodynamic regime, the particle distribution function relaxes to the Maxwellian 
	and the limiting system for particle density $n$ and particle macroscopic velocity $u$, which coincides to the fluid velocity, looks like a variable density incompressible Navier-Stokes system, see \eqref{eq:nu}. This justifies the asymptotic-preserving property of the scheme, which has a much relaxed numerical stability condition than a non-AP discretization, and can capture the hydrodynamic limit with time step and mesh size {\it independent} of the Stokes number.
	
	The paper is organized as follows.
	In Section \ref{sec:Model}, we detail some basic facts about the PDE system of interest to us, including its hydrodynamic limit system ($\varepsilon \rightarrow 0$). In Section \ref{sec:multiAP}, we give the details of the numerical scheme. Both first-order and second-order schemes are presented in this framework and the AP property will be studied. Section \ref{sec:numerical} is devoted to numerical simulation for checking accuracy, asymptotic behavior, and some applications, followed by conclusions in Section \ref{sec:conclusion}.

	\section{A Model Problem}
	\label{sec:Model}
	
	In this paper, we focus on the fine particle regime, in which
	the suitably scaled PDE systems for the multi-phase model  are given by \cite{JinLin2022}: 
	\begin{equation}\label{ModelEquation0}
		\left\{\begin{aligned}
			&\begin{aligned}
				(f_i)_{t}+v \cdot \nabla_{x} (f_i) - \nabla_{x} \Phi \cdot \nabla_{v} f_i &=\dfrac{1}{\varepsilon}\dfrac{1}{i^{2/3}}\mathcal{L}_{u,i} f_i, (t, x, v) \in \mathbb{R}^{+} \times \mathbb{T}^{3} \times \mathbb{R}^{3}, i=1,2,...,N, 
			\end{aligned}
			\\
			&u_{t}+\nabla_{x} \cdot ( u \otimes u)+\nabla_{x} p-\dfrac{1}{Re}\Delta_{x} u=\dfrac{\kappa}{\varepsilon} \sum_{i=1}^N\int_{\mathbb{R}^{3}}(v-u) f_i i^{1/3} \mbox{d} v, \quad(t, x) \in \mathbb{R}^{+} \times \mathbb{T}^{3}, \\
			&\nabla_{x} \cdot u=0,
		\end{aligned}\right.
	\end{equation}
	with the initial condition
	$$\left.u\right|_{t=0}=u_{0}, \quad \nabla_{x} \cdot u_{0}=0,\left.\quad f_i\right|_{t=0}=f_{i,0},$$
	where $\mathcal{L}_{u,i} f_i$ is the $i$th Fokker-Planck (FP) operator
	$$
	\mathcal{L}_{u,i} f_i=\nabla_{v} \cdot\left((v-u) f_{i}+\dfrac{\bar{\theta}}{i} \nabla_{v} f_{i}\right),\  i=1,2,\ldots,N,
	$$
	with $N$ the number of particle sizes.
	Without loss of generality, assume the reference temperature $\bar{\theta}=1$ throughout the paper. The discussion of the scaling issues is detailed in Appendix.
	Here $t\geq 0$ is  time, $x=\left(x_{1}, x_{2}\right) \in \Omega \subset \mathbb{R}^{2}$ is the space variable, and $v=\left(v_{1} , v_ {2}\right) \in \mathbb{R}^{2}$ is the particle velocity.
	$f_i = f_i(t,x,v), i=1,2,\ldots,N$ are the density functions of the particles. 
	$u=u(t, x)=\left(u_1(t, x), u_2(t, x)\right)$ is the velocity field of the fluid. 
	$\Phi=\Phi(x)$ is an external force field and  $\nabla_{x} \Phi$ represents the effect of the external force field on the particles. 
	The first equation describes the motion of particles. The two terms in the Fokker-Planck (FP) operator come from the drag force from the fluid and the effect of Brownian motions, respectively.  The second and third equations are the standard incompressible Navier-Stokes equations for the fluid, with the right-hand-side term describing the force coming from the particles.
	$\kappa>0$ is the coupling constant, which equals the ratio between the particle density $\rho_P$ and fluid density $\rho_F$,
	and $Re$ is the Reynolds number.
	$\varepsilon$ $(0<\varepsilon\leq 1)$ is the Stokes number given by $\varepsilon=\frac{2 \rho_p i^{2}}{9 \mu}$, with $\mu$ the dynamic viscosity of the fluid, $i$ the typical radius of the particles and $\rho_{P}$ the density of the particles. Here $\varepsilon$ is a constant with $i=1$. Its correlation with the particle size $i$ is given explicitly in the Eq. \ref{ModelEquation0}. 
	$\varepsilon=O(1)$ corresponds to the kinetic regime, while $\varepsilon\rightarrow 0$ corresponds to the fluid regime. 
	$Re$ is the Renalds number. 

	Let us briefly recall some basic facts about system \eqref{ModelEquation0} and the regime $\varepsilon \rightarrow 0$. The key remark, observed in \cite{GoudonJabin2004a, GoudonJabin2004b}, is the following energy-entropy dissipation property
	\begin{equation}\label{Energy}
		\begin{aligned}
			\frac{\mathrm{d}}{\mathrm{d} t}\left(\kappa\sum_{i=1}^{N} \int_{\mathbb{T}^{3}} \int_{\mathbb{R}^{3}} f_i\left(\ln (f_i)+1+i\Phi+i\frac{|v|^{2}}{2}\right) \,\mathrm{d} v \mathrm{d} x+\int_{\mathbb{T}^{3}}\frac{|\tilde{u}|^{2}}{2}\, \mathrm{d} x\right)+\int_{\mathbb{T}^{3}}\left|\nabla_{x} \tilde{u}\right|^{2} \,\mathrm{d} x\\
			+\frac{\kappa}{\varepsilon}\sum_{i=1}^{N} \int_{\mathbb{T}^{3}} \int_{\mathbb{R}^{3}}\left|(v-\tilde{u}) \sqrt{i^{1/3}f_i}+\frac{\bar{\theta} \nabla_{v} f_i}{\sqrt{i^{5/3}f_i}}\right|^{2} \, \mathrm{d} v \mathrm{d} x=0.
		\end{aligned}
	\end{equation}
	A similar relation holds when the problem is set on a bounded smooth domain $\Omega$ with reasonable boundary conditions. For instance one can assume no-slip of the fluid
	$$
	\left.u\right|_{\partial \Omega}=0
	$$
	and specular reflection of the particles
	$$
	\gamma^{-} f_i(t, x, v)=\gamma^{+} f_i(t, x, v-2 v \cdot \hat{v}(x) \hat{v}(x)),
	$$
	where $\hat{v}(x)$ stands for the unit outer normal at point $x \in \partial \Omega$ and $\gamma^{\pm}$denote the trace operators on the set
	$$
	\left\{(t, x, v) \in(0, \infty) \times \partial \Omega \times \mathbb{R}^{2}, \quad \pm v \cdot \hat{v}(x)>0\right\} .
	$$
	We refer to further comments in \cite{Carrilo2006}. It is worth rewriting the Fokker-Planck operator as
	$$
	L_{u,i} f_i=\frac{1}{i}\nabla_{v} \cdot\left(M_{u,i} \nabla_{v}\left(\frac{f_i}{M_{u,i}}\right)\right), \quad M_{i,u}(v)=\frac{i}{2 \pi} \exp \left(-\frac{i|v-u(t, x)|^{2}}{2}\right) .
	$$
	As $\varepsilon$ goes to $0$ , since the Fokker-Planck operator is penalized and the change of particles sizes is ignored, we expect that $f_i$ makes $L_{u,i} f_i$ (and the dissipation term in \eqref{Energy}) vanish which means that $f_i$ becomes proportional to the Maxwellian centered on the fluid velocity
	$$
	f_i(t, x, v) \simeq n_i(t, x) M_{i,u(t, x)} .
	$$
	Hence the question is to identify the equation satisfied as $\varepsilon \rightarrow 0$ by the particles density $n$ and the velocity $u$.
	
	For the deterministic multi-size particle-fluid systems \eqref{ModelEquation0}, we associate to $f_i(t,x,v), i=1,2,\ldots,N$ the following macroscopic quantities:
	$$
	\begin{gathered}
		n_{i}(t, x)=\int_{\mathbb{R}^{3}} f_i(t, x, v) \mathrm{d} v, \quad \rho_i(t,x) = i n_i(t,x), 
		\quad J_{i}(t, x)=  i \int_{\mathbb{R}^{3}} v f_i(t, x, v) \mathrm{d} v, \\
		\mathbb{P}_{i}(t, x)= i \int_{\mathbb{R}^{3}} v \otimes v f_i(t, x, v) \mathrm{d} v ,
	\end{gathered}
	$$
	where $\rho_i, J_i$ and $\mathbb{P}_{i}$ are the mass, momentum and stress tensors, respectively, of particles of size $i$.
	Integrating the first equation in \eqref{ModelEquation0} with respect to $i \mathrm{~d} v$ and $i v \mathrm{~d} v$ respectively, one obtains
	\begin{equation}
		i \partial_{t} n_{i}+\nabla_{x} \cdot J_{i}=0,
	\end{equation}
	and
	\begin{equation}\label{macroJ}
		\partial_{t}( J_{i})+\operatorname{Div}_{x} (\mathbb{P}_{i}) +i n_{i} \nabla_x \Phi =-\frac{1}{i^{2 / 3}\varepsilon} \left(J_{i}- i n_{i} u\right).
	\end{equation}
	Combined to system \eqref{ModelEquation0}, one has
	\begin{equation}\label{eq:kappaJ}
		\begin{array}{r}
			\partial_{t}\left(u+ \kappa \displaystyle\sum_{i=1}^{N} ( J_{i})\right)+\operatorname{Div}_{x}\left(u \otimes u+\kappa \displaystyle\sum_{i=1}^{N}( \mathbb{P}_{i})\right) 
			+ \nabla_{x} p \\
			+ \kappa \displaystyle\sum_{i=1}^{N} (i n_{i}) \nabla_x \Phi 
			-\dfrac{1}{Re}\Delta_{x} u=0.
		\end{array}	
	\end{equation}
	Accordingly, for $\varepsilon<<1$,  $J_i$ and $\mathbb{P}_i$ are asymptotically defined by the moments of the Maxwellian, i.e., 
	$$J_i \simeq i n_i u, \quad \mathbb{P}_i \simeq i n_i u \otimes u+ i n_i  \mathbb{I}.$$
	Inserting this ansatz into \eqref{eq:kappaJ}, one arrives at
	\begin{equation}
		\begin{aligned}
			\partial_{t}\left(\left(1+ \kappa \sum_{i=1}^{N}  i n_i\right)u\right)+\operatorname{Div}_{x}\left(\left(1+ \kappa \sum_{i=1}^{N}  i n_i\right)u \otimes u\right) + \nabla_{x} \left(p+ \kappa \sum_{i=1}^{N}  i n_i\right) &\\
			+\kappa \sum_{i=1}^{N} (i n_{i}) \nabla_x \Phi-\frac{1}{Re}\Delta_{x} u=0,&
		\end{aligned}	 
	\end{equation}
	Denote $\nu = \displaystyle\sum_{i=1}^{N}  i n_i$ with $n_{i}(t, x)=\int_{\mathbb{R}^{3}} f_i(t, x, v) \mathrm{d} v$. As $\varepsilon\rightarrow 0$,  \eqref{ModelEquation0} has a hydrodynamic limit
	\begin{equation}\label{eq:nu}
		\left\{\begin{aligned}
			& \partial_{t} \nu +\nabla_{x} \cdot(\nu u)=0, \\
			&\partial_{t}\left(\left(1+ \kappa \nu\right)u\right)+\operatorname{Div}_{x}\left(\left(1+ \kappa \nu\right)u \otimes u\right)+ \nabla_{x} \left(p+ \kappa \nu \right) + \kappa  \nu \nabla_x \Phi -\frac{1}{Re}\Delta_{x} u=0, \\
			&\nabla_{x} \cdot u=0,
		\end{aligned}\right.
	\end{equation}
	which is an incompressible Navier-Stokes system for the composite and inhomogeneous density $(1+\kappa \nu)$.
	
	In this paper, we are interested in numerical approximations of system \eqref{ModelEquation0}. We will pay particular attention to  the scaling parameter $\varepsilon$: the scheme should work over a wide range of values of the parameter, capturing the expected asymptotic behavior without introducing restrictions that would make small $\varepsilon$'s simulations numerically prohibitive. 
	This scheme consists in discretizing implicitly the stiff terms in the equations,  but it should be done as simple as possible because the inversion of the corresponding discrete system will be the main source of numerical cost.
	
	The scheme we developed can automatically capture the hydrodynamic limit \eqref{eq:nu} as $\varepsilon\rightarrow 0$. This is the so-called Asymptotic Preserving (AP) property, a term first introduced by Jin \cite{Jin1999}. The AP scheme is effective in the hydrodynamic regime ($\varepsilon \ll 1$) because it allows capturing the hydrodynamic limit \eqref{eq:nu} without numerically resolving the small scale $\varepsilon$. We refer to \cite{Jin2010, Jin2022} for reviews on AP schemes and their applications.

	\section{An AP Scheme for Multi-phase Flows}
	\label{sec:multiAP}
	
	We now give the details to update the numerical unknowns, having at hand $u^k, p^k, f_i^k$ and thus 
	$$n_i^k=\int_{\mathbb{R}^{d}} f_i^k \mathrm{d} v, \quad J_i^k=i \int_{\mathbb{R}^{d}} v f_i^k \mathrm{d} v,\quad i=1,2,...,N.
	$$
	
	\subsection{The time discretization - first order}
	
	\textbf{Step 1. Advancing densities. Compute $n_i^{k+1}$.}
	\begin{equation}\label{step1n}
		\begin{aligned}
			\dfrac{1}{\Delta t}\left(n_i^{k+1}-n_i^{k}\right)&=-\int v \cdot \nabla_{x} f_i^{k} \mathrm{~d} v
		\end{aligned}	
	\end{equation}
	
	\textbf{Step 2. Pressureless step. Compute $u^*, J_i^*$.}
	
	Notice that in order to derive an AP scheme, one has to impose the stiff drag force term in both the viscosity step and the projection step.
		
	Solve the viscosity part of momentum equations with only part of the stiff term:
	\begin{equation}\label{step2Ju}
		\begin{aligned}
			&\dfrac{1}{\Delta t}\left(J_i^{*}-J_i^{k}\right)=-i \int v \otimes v \nabla_{x} f_i^{k} \mathrm{~d} v-i n_i^{k} \nabla_{x} \Phi-\dfrac{1-\alpha}{i^{2/3}\varepsilon}\left(J_i^{*}-i n_i^{k+1} u^{*}\right), \\
			&\dfrac{1}{\Delta t}\left(u^{*}-u^{k}\right)-\dfrac{1}{Re}\Delta_{x} u^{*}=-\nabla_{x} \cdot\left(u^{k} \otimes u^{k}\right) +\sum_{i=1}^{N} \dfrac{(1-\alpha)\rho_P}{i^{2/3}\varepsilon\rho_F} \left(J_i^{*}-i n_i^{k+1} u^{*}\right),
		\end{aligned}
	\end{equation}
	where $\alpha \in(0,1)$ is any constant. One can simply choose $\alpha=\dfrac{1}{2}$.
	
	Eliminating $J_i^*$, and with the no-slip boundary condition for $u^{*}$ used, one obtains a Helmholtz equation for $u^*$:
	\begin{equation}\label{step2equ}
		\left\{ \begin{array}{l}
			\begin{aligned}
				&\left(\dfrac{1}{\Delta t}+\sum_{i=1}^{N} \dfrac{(1-\alpha)\rho_P}{(i^{2/3}\varepsilon+(1-\alpha) \Delta t)\rho_F} i n_i^{k+1}-\dfrac{1}{Re}\Delta_{x}\right) u^{*}\\&=\dfrac{u^{k}}{\Delta t}-\nabla_{x} \cdot\left(u^{k} \otimes u^{k}\right) + \sum_{i=1}^{N}\dfrac{(1-\alpha)\rho_P}{(i^{2/3}\varepsilon+(1-\alpha) \Delta t)\rho_F}\left(J_i^{k}-i \Delta t \int v \otimes v \nabla_{x} f_i^{k} \mathrm{~d} v-i \Delta t n_i^{k} \nabla_{x} \Phi\right)
			\end{aligned}
			\\
			\left.u^{*}\right|_{\partial \Omega}=0 
		\end{array}	\right.
	\end{equation}
	
	$u^*$ can be solved by the Preconditioned Counugate Gradient method and then $J_i^*$ can be solved accordingly from \eqref{step2Ju}.

	\textbf{Step 3. Projection step. Compute $p^{k+1}, u^{k+1}$.}
	
	Next $u^*$ is projected to the divergence free space, with the remaining stiff coupling term:
	\begin{equation}\label{step3equ0}
		\begin{aligned}
			&\dfrac{1}{\Delta t}\left(J_i^{* *}-J_i^{*}\right)=-\dfrac{\alpha}{i^{2/3}\varepsilon}\left(J_i^{* *}-i n_i^{k+1} u^{k+1}\right) \\
			&\dfrac{1}{\Delta t}\left(u^{k+1}-u^{*}\right)+\nabla_{x} p^{k+1}=\sum_{i=1}^{N} \dfrac{\alpha\rho_P}{i^{2/3}\varepsilon\rho_F}\left(J_i^{* *}-i n_i^{k+1} u^{k+1}\right)
		\end{aligned}
	\end{equation}
	i.e.,
	\begin{equation}\label{step3equ}
		\left(1+\sum_{i=1}^{N}\dfrac{\dfrac{\alpha\rho_P}{i^{2/3}\varepsilon\rho_F}}{\dfrac{1}{\Delta t}+\dfrac{\alpha}{i^{2/3}\varepsilon}}i n_i^{k+1} \right)u^{k+1}+ \Delta t \nabla_{x} p^{k+1}= u^{*} + \sum_{i=1}^{N} \dfrac{\dfrac{\alpha\rho_P}{i^{2/3}\varepsilon\rho_F}}{\dfrac{1}{\Delta t}+\dfrac{\alpha}{i^{2/3}\varepsilon}} J_i^{*}.
	\end{equation}
	
	Noting that $u^{k+1}$ is divergence free, by taking the divergence of both sides, one has
	\begin{equation}\label{step3p}
		\left\{ \begin{array}{l}
			\nabla_{x} \cdot\left(\dfrac{1}{\rho_\varepsilon^{k+1}} \nabla_{x} p^{k+1}\right)=\dfrac{1}{\Delta t} \nabla_{x} \cdot\left(\dfrac{u^{*} + \displaystyle\sum_{i=1}^{N} \dfrac{\dfrac{\alpha\rho_P}{i^{2/3}\varepsilon\rho_F}}{\dfrac{1}{\Delta t}+\dfrac{\alpha}{i^{2/3}\varepsilon}} J_i^{*}}{\rho_\varepsilon^{k+1}}\right),\\
			\left.\quad \dfrac{\partial p^{k+1}}{\partial \hat{v}}\right|_{\partial \Omega}=0,
		\end{array}	\right.
	\end{equation}
	with 
	$\rho_\varepsilon^{k+1} = 1+\displaystyle\sum_{i=1}^{N}\dfrac{\dfrac{\alpha\rho_P}{i^{2/3}\varepsilon\rho_F}}{\dfrac{1}{\Delta t}+\dfrac{\alpha}{i^{2/3}\varepsilon}}i n_i^{k+1}.$
	
	$p^{k+1}$ can be solved by a Conjugate Gradient method and then $u^{k+1}$ can be obatined accordingly by \eqref{step3equ0}.
	
	\textbf{Step 4. Kinetic equation. Compute $f_i^{k+1}, J_i^{k+1}$.}
	
	$f_i^{k+1}$ is solved based on the equation
	\begin{equation}\label{step4eqf}
		\begin{aligned}
			&\dfrac{f_i^{k+1}-f_i^{k}}{\Delta t}+v \cdot \nabla_{x} f_i^{k}-\nabla_{x} \Phi \cdot \nabla_{v} f_i^{k}=\dfrac{1}{i^{2/3}\varepsilon} \mathcal{L}_{u^{k+1},i} f_i^{k+1},
		\end{aligned}
	\end{equation}
	where
	$$
	\mathcal{L}_{u^{k+1},i} f_i^{k+1}=\frac{1}{i}\nabla_{v} \cdot\left(\left(v-u^{k+1}\right) f_i^{k+1}+\dfrac{1}{i}\nabla_{v} f_i^{k+1}\right) ,\quad i=1, 2,\cdots,N.
	$$
	Then $J_i^{k+1} (i=1,2,\cdots,N)$ is updated by taking the first moment of $f_i^{k+1} (i=1,2,\cdots,N)$.

	\subsection{The time discretization - second order}
	\label{Sec3.2}
	
	Now we generalize  the first order scheme \eqref{step1n}-\eqref{step4eqf} to second order. The convergence order can be improved by the following techniques.
	\begin{itemize}
		\item The time derivative terms are approximated by a second order BDF method, i.e.,
		$$
		\partial_{t} a\left(t^{k+1}\right) \approx \dfrac{3 a^{k+1}-4 a^{k}+a^{k-1}}{2 \Delta t}
		$$
		\item The transport terms are approximated by extrapolation from previous two steps, i.e.,
		$$
		b\left(t^{k+1}\right) \approx 2 b^{k}-b^{k-1} ;
		$$
		\item The stiff terms are implicitly evaluated at $t^{k+1}$;
		\item Take $\alpha=O(\Delta t)$ in the splitting of stiff terms.
		\item Pressure incremental technique: The viscosity step in momentum equation is solved with pressure at current step and then the projection step computes the pressure increment.
	\end{itemize}
	
	\textbf{Step 1. Advancing densities. Compute $n_i^{k+1}$.}
	\begin{equation}\label{step1n2}
		\begin{aligned}
			&\dfrac{1}{2 \Delta t}\left(3 n_i^{k+1}-4n_i^{k}+n_i^{k-1}\right)=-\int v \cdot \nabla_{x} (2f_i^{k}-f_i^{k-1}) \mathrm{~d} v.
		\end{aligned}	
	\end{equation}
	
	\textbf{Step 2. Pressureless step. Compute $u^*, J_i^*$.}
	
	Solve the viscosity part of momentum equations with only part of the stiff term:
	\begin{equation}\label{Step2J}
		\begin{aligned}
			&\frac{1}{2 \Delta t}\left(3 J_i^{*}-4 J_i^{k}+J_i^{k-1}\right)\\
			& \quad =-\int v \otimes v \nabla_{x} (2f_i^{k}-f_i^{k-1}) \mathrm{d} v-(2n_i^{k}-n_i^{k-1}) \nabla_{x} \Phi-\frac{1-\alpha}{i^{2/3}\varepsilon}\left(J_i^{*}-i n_i^{k+1} u^{*}\right) \\
			&\frac{1}{2 \Delta t}\left(3 u^{*}-4 u^{k}+u^{k-1}\right)-\dfrac{1}{Re}\Delta_{x} u^{*}+\nabla_{x} p^{k}\\
			&\quad = -\nabla_{x} \cdot (2(u^k \otimes u^k)-(u^{k-1} \otimes u^{k-1})) +\sum_{i=1}^{N}\frac{(1-\alpha)\rho_P}{i^{2/3}\varepsilon\rho_F} \left(J_i^{*}-i n_i^{k+1} u^{*}\right)
		\end{aligned}
	\end{equation}
	Again $\alpha \in(0,1) .$ We need $\alpha=O(\Delta t)$ to ensure the second order accuracy.

	Eliminating $J_i^*$, and with the no-slip boundary condition for $u^{*}$ used, one obtains a Helmholtz equation for $u^*$:
	\begin{equation}\label{step2equ2}
		\left\{ \begin{array}{l}
			\begin{aligned}
				&\left(\frac{3}{2 \Delta t}+ \sum_{i=1}^{N}\frac{3(1-\alpha)\rho_P}{(3 *i^{2/3}\varepsilon+2(1-\alpha) \Delta t)\rho_F} i n_i^{k+1} -\dfrac{1}{Re}\Delta_{x}\right) u^{*} \\
				&\quad =  \frac{4 u^{k}-u^{k-1}}{2 \Delta t}-\nabla_{x} \cdot(2(u^k \otimes u^k)-(u^{k-1} \otimes u^{k-1})) -\nabla_{x} p^{k}\\
				&\quad\quad +\sum_{i=1}^{N}\frac{(1-\alpha)\rho_P}{(3*i^{2/3}\varepsilon+2(1-\alpha) \Delta t)\rho_F} \bigl\{4 J_i^{k}-J_i^{k-1} \bigr.\\
				& \quad\quad\quad\quad \bigl. -2i \Delta t\int v \otimes v \nabla_{x} (2f_i^{k}-f_i^{k-1}) \mathrm{d} v -2i \Delta t (2n_i^{k}-n_i^{k-1}) \nabla_{x} \Phi\bigr\}.
			\end{aligned}
			\\
			\left.u^{*}\right|_{\partial \Omega}=0 
		\end{array}	\right.
	\end{equation}
	
	$u^*$ can be solved by the Preconditioned Counugate Gradient method and then $J^*$ can be solved accordingly.
	\begin{equation*}
		\begin{aligned}
			J_i^{*} =& \frac{i^{2/3}\varepsilon}{3*i^{2/3}\varepsilon +(1-\alpha)2\Delta t }\left[4 J_1^{k}-J_1^{k-1}-2i\Delta t\int v \otimes v \nabla_{x} (2f_i^{k}-f_i^{k-1}) \mathrm{d} v\right.\\
			&\left.-2i\Delta t(2n_i^{k}-n_i^{k-1}) \nabla_{x} \Phi+\frac{(1-\alpha)2\Delta t}{i^{2/3}\varepsilon}i n_i^{k+1} u^{*} \right].
		\end{aligned}
	\end{equation*}
	
	\textbf{Step 3. Projection step. Compute $p^{k+1}, u^{k+1}$.}
	
	Next $u^*$ is projected to the divergence free space, with the remaining stiff coupling term:
	\begin{equation}\label{step3J}
		\begin{aligned}
			&\frac{3}{2 \Delta t}\left(J_i^{* *}-J_i^{*}\right)=-\frac{\alpha}{i^{2/3}\varepsilon}\left(J_i^{* *}-i n_i^{k+1} u^{k+1}\right), \\
			&\frac{3}{2 \Delta t}\left(u^{k+1}-u^{*}\right)+\nabla_{x}\left(p^{k+1}-p^{k}\right)= \sum_{i=1}^{N}\frac{\alpha\rho_P}{i^{2/3}\varepsilon\rho_F} \left(J_i^{* *}-i n_i^{k+1} u^{k+1}\right) .
		\end{aligned}
	\end{equation}
	i.e.,
	\begin{equation}\label{step3JJ}
		\left(1 +\sum_{i=1}^{N} \dfrac{\frac{\alpha\rho_P}{i^{2/3}\varepsilon\rho_F}}{\frac{3}{2 \Delta t}+\frac{\alpha}{i^{2/3}\varepsilon}}i n_i^{k+1} \right)u^{k+1}+ \frac{2 \Delta t}{3} \nabla_{x}\left(p^{k+1}-p^{k}\right)= u^{*}+\sum_{i=1}^{N} \dfrac{\frac{\alpha\rho_P}{i^{2/3}\varepsilon\rho_F}}{\frac{3}{2 \Delta t}+\frac{\alpha}{i^{2/3}\varepsilon}}  J_i^{*}  .
	\end{equation}
	
	Noting that $u^{k+1}$ is divergence free, by taking the divergence of both sides,
	\begin{equation}\label{step3p2}
		\left\{ \begin{array}{l}
			\nabla_{x} \cdot\left(\dfrac{1}{\rho_\varepsilon^{k+1}} \nabla_{x} (p^{k+1}-p^k)\right)=\dfrac{3}{2\Delta t} \nabla_{x} \cdot\left(\dfrac{u^{*}+\displaystyle\sum_{i=1}^{N} \dfrac{\frac{\alpha\rho_P}{i^{2/3}\varepsilon\rho_F}}{\frac{3}{2 \Delta t}+\frac{\alpha}{i^{2/3}\varepsilon}}  J_i^{*}}{\rho_\varepsilon^{k+1}}\right),\\
			\left.\quad \dfrac{\partial p^{k+1}}{\partial \hat{v}}\right|_{\partial \Omega}=0,
		\end{array}	\right. 
	\end{equation}
	with $\rho_\varepsilon^{k+1} = 1 +\displaystyle\sum_{i=1}^{N} \dfrac{\frac{\alpha\rho_P}{i^{2/3}\varepsilon\rho_F}}{\frac{3}{2 \Delta t}+\frac{\alpha}{i^{2/3}\varepsilon}}i n_i^{k+1}$.
	
	$p^{k+1}$ can be solved by a Conjugate Gradient method and then $u^{k+1}$ can be obatined accordingly.
	
	\textbf{Step 4. Kinetic equation. Compute $f_i^{k+1}, J_i^{k+1}$.}
	
	$f_i^{k+1}$ is solved based on the equation
	\begin{equation}\label{step4f2}
		\begin{aligned}
			&\frac{3 f_i^{k+1}-4 f_i^{k}+f_i^{k-1}}{2 \Delta t}+\left(v \cdot \nabla_{x}-\nabla_{x} \Phi \cdot \nabla_{v}\right)\left(2 f_i^{k}-f_i^{k-1}\right)=\frac{1}{i^{2/3}\varepsilon} \mathcal{L}_{u^{k+1},i} f_i^{k+1}
		\end{aligned}	
	\end{equation}
	where
	$$
	\mathcal{L}_{u^{k+1},i} f_i^{k+1}=\frac{1}{i}\nabla_{v} \cdot\left(\left(v-u^{k+1}\right) f_i^{k+1}+\dfrac{1}{i}\nabla_{v} f_i^{k+1}\right),\quad i=1, 2,\cdots,N.
	$$
	Then $J_i^{k+1} (i=1, 2,...,N)$ is updated by taking the first moment of $f_i^{k+1} (i=1, 2,...,N)$.
	
	Eqs. \eqref{step1n2}-\eqref{step4f2} give a second order scheme in time. We will check this convergence order numerically in Section \ref{sec:conver}. Note that this second order scheme is a multi-step method. To compute the solutions at $t^{k+1}$, we need the solutions from both $t^k$ and $t^{k-1}$. Therefore, with initial data at $t^0$, it is necessary to apply a first order method to obtain the solutions at $t^1$. This second order scheme can then be started.

	\subsection{Space and velocity discretizations}
	For the sake of completeness, let us discuss space and velocity discretizations by restricting to the two-dimension case. The extension to higher dimension is straightforward. Only Cartesian grids are considered. We denote by $\Delta x$ the (uniform) mesh size. We define a regularly spaced and symmetric velocity grid, with step $\Delta v$. Denoting $\mathbf{j}=\left(j, j^{\prime}\right)$ and $\mathbf{m}=\left(m, m^{\prime}\right)$ in $\mathbb{N}^2, f_{\mathbf{j}, \mathbf{m}}^k$ stands for the numerical approximation of $f\left(k \Delta t,(\mathbf{j}-\mathbf{1} / \mathbf{2}) \Delta x,(\mathbf{m}-\mathbf{1} / \mathbf{2}) \Delta v-v_{\max }\right)$. Here we assume $v \in\left[-v_{\max } v_{\max }\right]^2$. The grid points are located in the cell center.
	
	For the boundary condition, the specular reflection law is used to define the ghost points. For instance, labeling the numerical unknown with indices $j, j^{\prime} \in\{1, \ldots, J\}$ and $m, m^{\prime} \in\{1, \ldots, 2 M\}$, leads to
	$$
	f_{0, j^{\prime} ; m, m^{\prime}}^k=f_{1, j^{\prime} ; 2 M+1-m ; m^{\prime}}^k \quad f_{J+1, j^{\prime} ; m, m^{\prime}}^k=f_{J, j^{\prime} ; 2 M+1-m, m^{\prime}}^k .
	$$
	For the pressure, Neumann boundary condition in \eqref{step3p} and \eqref{step3p2} leads to
	$$
	p_{0, j^{\prime}}^{k+1}=p_{1, j^{\prime}}^{k+1}, \quad p_{J+1, j^{\prime}}^{k+1}=p_{J, \jmath^{\prime}}^{k+1} .
	$$
	The no-slip boundary of $u^*$ in \eqref{step2equ} and \eqref{step2equ2} leads to
	$$
	u_{0, j^{\prime}}^*=-u_{1, j^{\prime}}^*, \quad u_{J+1, j^{\prime}}^*=-u_{J, j^*}^*
	$$
	Similar expression holds when exchanging the role of $u_1^*, u_2^*, j, j$ and $m^{\prime}, m^{\prime}$.
	We refer for instance to \cite{Aregba2004} for discussion of numerical boundary conditions for kinetic schemes.	
	
	For the transport term $v \cdot \nabla_x f_i$ in \eqref{step4eqf} and \eqref{step4f2} and the derivative with respect to velocity which appears in the acceleration term, the upwind type second order shock capturing schemes with a kind of slope limiter is applied (see \cite{Goudon2012}). Discrete differential operators in higher dimension are defined dimension-by-dimension.	
	The convection term $\nabla_x \cdot(u \otimes u)$ and the diffusion term $\Delta_x u$ in incompressible Navier-Stokes system \eqref{ModelEquation0}, as well as the terms $\nabla_x p$ and $\nabla_x \cdot u^*$  in the projection steps \eqref{step3JJ} - \eqref{step3p2}, are approximated by centered differences.	
	Macroscopic quantities are defined by using the two-dimensional version of the trapezoidal rule in order to ensure that the even moments of the odd functions with respect to $v$ vanish.	
	
	The numerical stability analysis of the complete problem is beyond the scope of this paper. However, we can expect, and as confirmed by
	our numerical observations, that the only constraint on the time step is the CFL condition coming from the transport part of the kinetic equation \eqref{ModelEquation0}, i.e. $\Delta t \leqslant \frac{\Delta x}{\max |v|} $, with $\Delta x $ the space mesh size. 
	
	\subsection{Treatment of the Fokker-Planck Equation}
	
	Now we focus on how to solve $f_i^{k+1}$ from \eqref{step4eqf} and \eqref{step4f2} where the stiff term is treated implicitly.
	For each Fokker-Planck equation $$\partial_{t} f_i+v \cdot \nabla_{x} f_i=\dfrac{1}{i^{2/3}\varepsilon} \mathcal{L}_{u,i} f_i+ \nabla_{x} \Phi \cdot \nabla_{v} f_i, \quad i=1,2,\ldots,N,$$
	we introduce the "local Maxwellian":
	$$
	M_{u, i}(v):=\frac{i }{2 \pi}\exp \left(-\frac{i|v-u|^{2}}{2 }\right), \quad i=1,2,\ldots,N.
	$$
	The crucial observation consists in rewriting
	($\mathcal{L}_{u, i} f_i=\nabla_{v}\cdot \left((v-u) f_i+\frac{1}{i} \nabla_{v} f_i\right)$) 
	$$	\mathcal{L}_{u, i} f_i=\frac{1}{i} \nabla_{v}\cdot \left(M_{u, i} \nabla_{v}\left(\frac{f_i}{M_{u, i}}\right)\right), \quad i=1,2,\ldots,N. $$
	We need to invert the Fokker-Planck operators respectively. To this end, we follow the approach introduced in \cite{JinYan2011}.
	For the $i$-th $(1\leq i \leq N)$ Fokker-Planck equation,
	We write
	$$
	\mathcal{L}_{u,i}f_i=\sqrt{M_{u,i}} \tilde{\mathcal{L}}_u h_i
	$$
	with
	$$
	h_i=\frac{f_i}{\sqrt{M_{u,i}}}, \quad \tilde{\mathcal{L}}_{u,i} h_i=\frac{1}{i\sqrt{M_{u,i}}} \nabla_\nu \cdot\left(M_{u,i} \nabla_v\left(\frac{h_i}{\sqrt{M_{u,i}}}\right)\right) .
	$$
	Note that $\tilde{\mathcal{L}}_{u,i}$ is symmetric for the standard $L^2$ inner product
	$$
	\int_{\mathbb{R}^N} \tilde{\mathcal{L}}_u h g \mathrm{~d} v=\int_{\mathbb{R}^N} h \tilde{\mathcal{L}}_u g \mathrm{~d} v .
	$$

	We set
	$$
	h_{i,\mathbf{j} \mathbf{m}}=\frac{f_{i,\mathbf{j} \mathbf{m}}^{k+1}}{\sqrt{M_{i,\mathbf{j} \mathbf{m}}^{k+1}}}, \quad \mathscr{L} f_{i,\mathbf{j} \mathbf{m}}^{k+1}=
	\sqrt{M_{i,\mathbf{j} \mathbf{m}}^{k+1}} \tilde{\mathscr{L}}{h_{i,\mathbf{j} \mathbf{m}}}.
	$$
	The discrete operator $\mathscr{L}$ is symmetric which allows to make use of the Conjugate Gradient algorithm.
	In the two-dimension setting, the discrete operator $\tilde{\mathscr{L}}$ is defined as follows:
	\begin{equation}
		\begin{aligned}
			&\tilde{\mathscr{L}} h_{j, j^{\prime} ; m, m^{\prime}}=\frac{1}{\Delta v^{2}}\left(h_{j, j^{\prime} ; m, m^{\prime}+1}+h_{j, j^{\prime} ; m+1, m}-\bar{[M_i]}_{j, j^{\prime} ; m, m^{\prime}}^{k+1} h_{j, j^{\prime} ; m, m^{\prime}}+h_{j, j^{\prime} ; m, m^{\prime}-1}+h_{j, j^{\prime} ; m-1, m^{\prime}}\right), \\
			&\bar{[M_i]}_{j, j^{\prime} ; m, m^{\prime}}^{k+1}=\frac{\sqrt{[M_i]_{j, j^{\prime}, m+1, m^{\prime}}^{k+1}}+\sqrt{[M_i]_{j, j^{\prime} ; m, m^{\prime}+1}^{k+1}}+\sqrt{[M_i]_{j, j^{\prime} ; m-1, m^{\prime}}^{k+1}}+\sqrt{[M_i]_{j, j^{\prime} ; m, m^{\prime}-1}^{k+1}}}{\sqrt{[M_i]_{j, j^{\prime}, m, m^{\prime}}^{k+1}}}
		\end{aligned}
	\end{equation}
	which indeed leads to a symmetric matrix. Observe that $\tilde{\mathscr{L}}\left(\sqrt{M^{k+1}}\right) \quad=0$. 
	Therefore, the update of the particle distribution function consists of the following two steps:
	\begin{itemize}
		\item Step 1. Solve the linear system
		$$
		\left(1-\frac{\Delta t}{i^{5/3}\varepsilon}  \tilde{\mathscr{L}}\right) h_{i,\mathbf{j} \mathbf{m}}=\frac{f_{i,\mathbf{j} \mathbf{m}}^{k}-\Delta tv D_{x}[f_i]_{\mathbf{j} \mathbf{m}}^{k}+\Delta t v D_{x}[\Phi]_{\mathbf{j} \mathbf{m}}^{k} \frac{1}{v} D_{v}[f_i]^{k}_{\mathbf{j} \mathbf{m}}}{\sqrt{M_{i,\mathbf{j} \mathbf{m}}^{k+1}}},
		$$
		or
		$$
		\begin{aligned}
			&\left(1-\dfrac{2\Delta t}{3*i^{5/3}\varepsilon} \tilde{\mathscr{L}} \right) {h_i}_{\mathbf{j} \mathbf{m}}\\
			&=\dfrac{4f_{i,\mathbf{j} \mathbf{m}}^{k}-{f_i}_{\mathbf{j} \mathbf{m}}^{k-1}-2\Delta tv D_{x}(2[f_i]_{\mathbf{j} \mathbf{m}}^{k}-[f_i]_{\mathbf{j} \mathbf{m}}^{k-1})+2\Delta t v D_{x}[\Phi]_{\mathbf{j} \mathbf{m}}^{k} \frac{1}{v} D_{v}(2[f_i]_{\mathbf{j} \mathbf{m}}^{k}-[f_i]_{\mathbf{j} \mathbf{m}}^{k-1})}{3\sqrt{{M_i}_{\mathbf{j} \mathbf{m}}^{k+1}}}.
		\end{aligned}$$
		\item Step 2. Set $f_{i \mathbf{j} \mathbf{m}}^{k+1}=h_{i,\mathbf{j} \mathbf{m}} \sqrt{[M_i]_{\mathbf{j} \mathbf{m}}^{k+1}}$.
	\end{itemize}

	\subsection{Properties of the scheme}
	Now we show that the first order scheme \eqref{step1n}-\eqref{step4eqf} is asymptotic preserving and the limiting scheme gives a first order approximation for the limiting system \eqref{eq:nu}.
	As $\varepsilon \rightarrow 0$, \eqref{step4eqf} gives
	$$
	\mathcal{L}_{u^{k+1},i} f_i^{k+1}=O(\varepsilon), \quad \text { for } k \geqslant 0, \quad i=1,2,\ldots,N.
	$$
	This is equivalent to
	$$
	f_i^{k}=n_i^{k} M_{u^{k},i}+O(\varepsilon), \quad \text { for } k \geqslant 1,  \quad i=1,2,\ldots,N.
	$$
	Then one has
	$$
	\begin{aligned}
		&J_i^{k}=i n_i^{k} u^{k}+O(\varepsilon), \\
		&i \int_{\mathbb{R}^{N}} v \otimes v f_i^{k} \mathrm{~d} v=i n_i^{k} u^{k} \otimes u^{k}+i n_i^{k} \mathbb{I}+O(\varepsilon) .
	\end{aligned}
	$$
	Therefore, \eqref{step1n} is just
	\begin{equation}\label{APn}
		\frac{1}{\Delta t}\left(n_i^{k+1}-n_i^{k}\right)=-\nabla_{x} \cdot\left(n_i^{k} u^{k}\right)+O(\varepsilon) .
	\end{equation}
	Besides, the first and second Equation in \eqref{step2Ju} give
	$$
	J_i^{*}=i n_i^{k+1} u^{*}+O(\varepsilon) .
	$$
	Multiplying the first equation in \eqref{step2Ju} by $\kappa$, summing over $i$, and adding to the second equation in \eqref{step2Ju}, one obtains
	\begin{equation}\label{AP2}
		\begin{aligned}
			&\frac{1}{\Delta t}\left(\left(1+\kappa\sum_{i=1}^N i n_i^{k+1}\right) u^{*}-\left(1+\kappa \sum_{i=1}^N i n_i^{k}\right) u^{k}\right)-\frac{1}{Re}\Delta_{x} u^{*}\\
			&=-\nabla_{x} \cdot\left(\left(1+\kappa \sum_{i=1}^N i n_i^{k}\right) u^{k} \otimes u^{k}\right)-\kappa \nabla_{x} \sum_{i=1}^N i n_i^{k}-\kappa \sum_{i=1}^N i n_i^{k} \nabla_{x} \Phi+O(\varepsilon) .
		\end{aligned}
	\end{equation}
	Eqs. \eqref{APn} and \eqref{AP2} give a first order discretization of the limiting system \eqref{eq:nu}, without the pressure term. Moreover, as $\varepsilon \rightarrow 0$, \eqref{step3equ} becomes,
		$$
		u^{k+1}+\frac{1}{1+\kappa \sum_{i=1}^N i n_i^{k+1}} \Delta t \nabla_{x} p^{k+1}=u^{*},
		$$
	which is exactly the projection step for \eqref{eq:nu}.
	Similarly one can show the $\varepsilon \rightarrow 0$ limit of \eqref{step1n2}-\eqref{step4f2} is
	$$
	\begin{aligned}
		&\frac{1}{2 \Delta t}\left(3 n_i^{k+1}-4 n_i^{k}+n_i^{k-1}\right)=-\nabla_{x} \cdot\left(n_i^{k} u^{k}\right)^{\dagger} \\
		&\frac{1}{2 \Delta t}\left(3\left(1+\kappa \sum_{i=1}^N i n_i^{k+1}\right) u^{*}-4\left(1+\kappa \sum_{i=1}^N i n_i^{k}\right) u^{k}+\left(1+\kappa \sum_{i=1}^N i n_i^{k-1}\right) u^{k-1}\right)-\Delta_{x} u^{*}+\nabla_{x} p^{k}\\
		& \quad \quad =-\nabla_{x} \cdot\left(\left(1+\kappa \sum_{i=1}^N i n_i\right) u \otimes u\right)^{\dagger}-\kappa \nabla_{x} \sum_{i=1}^N i n_i^{\dagger}-\kappa \sum_{i=1}^N i n_i^{\dagger} \nabla_{x} \Phi, \\
		&\frac{3\left(u^{k+1}-u^{*}\right)}{2 \Delta t}+\frac{1}{1+\kappa \sum_{i=1}^N i n_i^{k+1}} \nabla_{x}\left(p^{k+1}-p^{k}\right)=0, \\
		&\nabla_{x} \cdot u^{k+1}=0 .
	\end{aligned}
	$$
	It is the second order projection scheme described in Section \ref{Sec3.2} for the limiting system \eqref{eq:nu}, an incompressible Navier-Stokes system with spatial variable density $(1+\kappa \nu)$. Here $\nu = \sum_{i=1}^N i n_i$ is defined as before. 
	Notice that it is consistent with the kinetic-fluid two-phase flow system provided that $N=1$ \cite{Goudon2013asymptotic}.
	
	\begin{remark} We can formally check the second order accuracy. 
		First, \eqref{Step2J} is (at least) a first order time discretization of the system \eqref{macroJ} and \eqref{ModelEquation0}. The local truncation error gives,
		\begin{equation}\label{reuJ}
			u^{*}=u^{k+1}+O\left(\Delta t^{2}\right), \quad J_i^{*}=J_i^{k+1}+O\left(\Delta t^{2}\right) .
		\end{equation}
		Next we add up Eqs. \ref{Step2J} and \eqref{step3J}:
		$$
		\begin{aligned}
			&\frac{1}{2 \Delta t}\left(3 J_i^{**}-4 J_i^{k}+J_i^{k-1}\right) \\
			&=-i\int v \otimes v \nabla_{x} f_i^{\dagger} \mathrm{d} v-i n_i^{\dagger} \nabla_{x} \Phi-\frac{1}{i^{2/3}\varepsilon}\left(J_i^{**}-i n_i^{k+1} u^{k+1}\right)+R_{1,i},\  i=1, 2, ..., N,\\
			&\frac{1}{2 \Delta t}\left(3 u^{k+1}-4 u^{k}+u^{k-1}\right)-\frac{1}{Re}\Delta_{x} u^{k+1}+\nabla_{x} p^{k+1} \\
			&=-\nabla_{x} \cdot(u \otimes u)^{\dagger}+\kappa\sum_{i=1}^N\frac{1}{i^{2/3}\varepsilon} \left(J_i^{**}-i n_i^{k+1} u^{k+1}\right)+R_{2},
		\end{aligned}
		$$
		where $b^\dagger = 2 b^k - b^{k-1}$, and the remainder terms are given by
		$$
		\begin{aligned}
			&R_{1,i}=-\frac{1-\alpha}{\varepsilon}\left(\left(J_i^{*}-i n_i^{k+1} u^{*}\right)-\left(J_i^{**}-i n_i^{k+1} u^{k+1}\right)\right), \\
			&R_{2}=\frac{1-\alpha}{\varepsilon} \kappa\left(\left(J^*_i-i n_i^{k+1} u^{*}\right)-\left(J_i^{**}-i n_i^{k+1} u^{k+1}\right)\right)+\Delta_{x}\left(u^{*}-u^{k+1}\right) .
		\end{aligned}
		$$
		Noting that \eqref{step3J} combined with \eqref{Step2J} is also (at least) a first order time discretization of the system \eqref{macroJ} and \eqref{ModelEquation0}, one has $J^{**}=J^{k+1}+O\left(\Delta t^{2}\right)$. Combined with \eqref{reuJ}, one has
		$$
		R_{1,i}=O\left(\Delta t^{2}\right), \quad R_{2}=O\left(\Delta t^{2}\right) .
		$$
		Therefore $u^{k+1}$ and $J_i^{* *}$ are second order approximations of $u\left(t^{k+1}\right)$ and $J_i\left(t^{k+1}\right)$. Then the distribution $f_i^{k+1}$ is solved via the second order discretization \eqref{step4f2}.
	\end{remark}

	\section{Numerical Simulation}
	\label{sec:numerical}
	
	Let us now check the performances of the method through a set of numerical experimetns.
	From now on, we will use the following notation: $\mathbf{x}=(x, y)$ is the position variable, $\mathbf{v}=\left(v_1, v_2\right)$ is the velocity variable, $\mathbf{u}=\left(u_1, u_2\right)$ is the fluid velocity, $\mathbf{u}_p=\left(u_{p 1}, u_{p 2}\right)$ is the macroscopic particle velocity. 
	
	We always use the following settings unless otherwise stated.

	\begin{itemize}
		
		\item The computation is performed on $(\mathbf{x}, \mathbf{v}) \in[0,1]^2 \times\left[-v_{\max }, v_{\max }\right]^2 $ , with $v_{\max }=8$. 
		We take $N_x=128$ grid points in each $x$ direction and $N_v=32$ grid points in each $v$ direction.
		Denoting $\mathbf{j}=\left(j, j^{\prime}\right)$ and $\mathbf{m}=\left(m, m^{\prime}\right)$ in $\mathbb{N}^2, f_{i,\mathbf{j} ; \mathbf{m}}^k$ stands for the numerical approximation of $f_i\left(k \Delta t,(\mathbf{j}-\mathbf{1} / \mathbf{2}) \Delta x,(\mathbf{m}-\mathbf{1} / \mathbf{2}) \Delta v-v_{\max }\right)$.

		\item  For the boundary condition, labeling the numerical unknown with indices $j, j^{\prime} \in\{1, \ldots, J\}$ and $m, m^{\prime} \in\{1, \ldots, 2 M\}$ where the $M$ first (resp. last) velocities are negative (resp. positive), the specular reflection for particle distributions $f$, the no-flip boundary condition for fluid velocity $\mathbf{u}$ and Neumann boundary condition for the pressure $p$ lead to
		$$
		f_{0, j^{\prime} ; m, m^{\prime}}^k=f_{1, j^{\prime} ; 2 M+1-m ; m^{\prime}}^k, \quad f_{J+1, j^{\prime} ; m, m^{\prime}}^k=f_{J, j^{\prime} ; 2 M+1-m, m^{\prime}}^k .
		$$
		$$
		u_{0, j^{\prime}}^*=-u_{1, j^{\prime}}^*, \quad u_{j+1, j^{\prime}}^*=-u_{j, f^*}^* \text {. }
		$$
		$$
		p_{0, j^{\prime}}^{k+1}=p_{1, J^{\prime}}^{k+1}, \quad p_{J+1, J^{\prime}}^{k+1}=p_{J, j^{\prime}}^{k+1} .
		$$
		Similar expression holds when exchanging the role of $u_1^*, u_2^*, j, j^{\prime}$ and $m, m^{\prime}$.

		\item We apply the second-order method described in Section \ref{Sec3.2}.  
		For the transport term $v \cdot \nabla_x f$ in \eqref{step4eqf} and  \eqref{step4f2}, the upwind type second order shock capturing scheme with a van Leer type slope limiter is applied (see \cite{Goudon2012}).
		The convection term $\nabla_x \cdot(u \otimes u)$ and the diffusion term $\Delta_x u$ in incompressible Navier-Stokes system \eqref{ModelEquation0}, as well as the terms $\nabla_x p$  and $\nabla_x \cdot u^*$ in the projection steps \eqref{step3equ}-\eqref{step3p}, are approximated by centered differences.
		Macroscopic quantities are defined by using the 2-dimensional version of the trapezoidal rule in order to ensure that the even moments of the odd functions with respect to $v$ vanish.
		For the derivative with respect to velocity which appears in the acceleration term, the upwind type second order shock capturing scheme is applied (see \cite{Goudon2012}).

		\item The time step is taken by $\Delta t=\frac{\Delta x}{5 v_{\max }}$ to ensure the numerical stability.

		\item We always take
		$$
		f_i(0, \mathbf{x}, \mathbf{v})=n_i(0, \mathbf{x}) M_{\mathbf{u}_p(0, \mathbf{x}),i}
		$$
		as initial data for particle distributions. Here $\mathbf{u}_{p,i}=\frac{J_i}{n_i}$ is the macroscopic velocity of the particle. We point out that this is not necessarily an equilibrium state when  $\mathbf{u}_p \neq \mathbf{u}$ and thus $\mathcal{L}_{\mathbf{u},i} f_i \neq 0$ in \eqref{ModelEquation0}.
		That is, for initial data, we do not require $\mathbf{u}_{p,0}=\mathbf{u}_0$.
		
		\item We take $\kappa=2$ throughout the simulations. Note that the scheme can be applied to the case with far larger values of $\kappa$ without any difficulty.
		
		\item If the gravity is taken into account, we take $\Phi=g y$ with gravity constant $g=1$.
		
		\item Our simulation allows us to compute with Reynolds number up to Re $=1000$ without trouble in stability. $Re=1$ and $Re=1000$ are both taken into account in simulations.
		
		\item Our AP scheme allows us to simulate the multi-phase partidulate flows with any finite number of distinct size of particles. For the sake of clarity, we take $N=2$ throughout the simulation, that is, we consider a three-phase system made up of a smaller particle ($i=1$), a bigger particle ($i=2$) and the fluid.
		
	\end{itemize}

	\subsection{Convergence order}
	\label{sec:conver}
	
	First we numerically check the accuracy of the scheme described in Section \ref{Sec3.2}. We start with the initial data
	\begin{equation}\label{data0}
		\begin{aligned}
			&n_1(0, \mathbf{x})=n_2(0, \mathbf{x})=10^{-10}+\exp \left(-80(x-0.5)^{2}-80(y-0.5)^{2}\right), \\
			&\mathbf{u}_{p,1}(0, \mathbf{x})=\mathbf{u}_{p,2}(0, \mathbf{x})=\left(\begin{array}{c}
				\sin ^{2}(\pi x) \sin (2 \pi y) \\
				-\sin ^{2}(\pi y) \sin (2 \pi x)
			\end{array}\right), \\
			&\mathbf{u}(0, \mathbf{x})=\mathbf{u}_{p,1}(0, \mathbf{x}) .
		\end{aligned}
	\end{equation}
	We compute the solutions on a grid of $N_{x} \times N_{x} \times N_{v} \times N_{v}$. We set $\Delta x = \frac{1}{N_x}$ with $N_{x}=16,32,64,128$ respectively, and $N_{v}=32$. At the final time $t_{\max }=0.025$, we check the following realtive error in $\ell^p$ norm,
	$$
	\begin{aligned}
		&e_{\Delta x}(f_i)=\max _{t \in\left(0, t_{\max }\right)} \frac{\left\|f_{i,\Delta x}(t)-f_{i,2 \Delta x}(t)\right\|_{p}}{\left\|f_{i,2 \Delta x}(0)\right\|_{p}}, \\
		&e_{\Delta x}(\mathbf{u})=\max _{t \in\left(0, t_{\max }\right)} \frac{\left\|\mathbf{u}_{\Delta x}(t)-\mathbf{u}_{2 \Delta x}(t)\right\|_{p}}{\left\|\mathbf{u}_{2 \Delta x}\left(t_{\max }\right)\right\|_{p}} .
	\end{aligned}
	$$
	This can be considered as an estimation of the relative error in $\ell^p$ norm, where $f_{i,\Delta x}$ and $\mathbf{u}_{\Delta x}$ are the numerical solutions computed from a grid of size $\Delta x=\frac{1}{N_{x}}$.
	We bear in mind that the stability constraint imposes $\Delta t=O(\Delta x)$.
	We shall say that the numerical scheme is of order $k$ if $e_{\Delta x} \leq C \Delta x^k$ holds,  for $\Delta x$ small enough. 
	Here simulations are performed with the Reynolds number $\mathrm{Re}=1$.

	\begin{figure}[htbp]
		\centering
		\includegraphics[width=15cm]{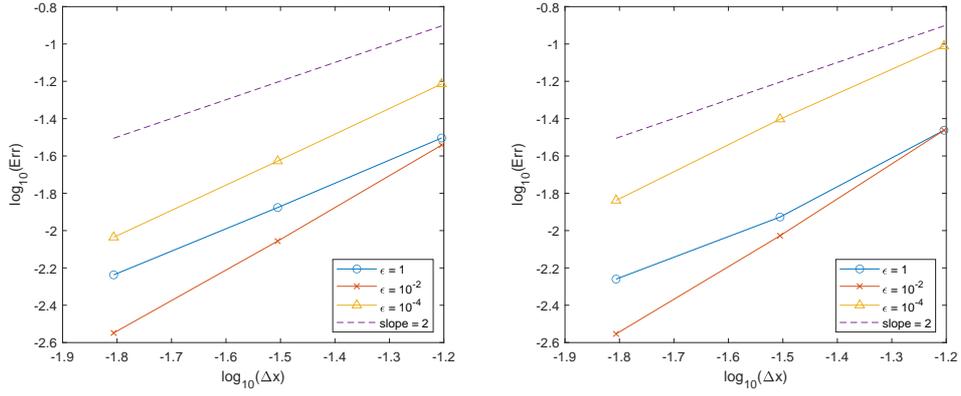}\\
		\caption{The test of convergence order with initial data \eqref{data0}. This figure shows the $l^2$ norm in the  small particle distribution $f_1$ (left) and the large particle distribution $f_2$ (right) with different $\varepsilon$.}\label{fig:ThreeEx0f}
	\end{figure}
	
	\begin{figure}[htbp]
		\centering
		\includegraphics[width=6.6cm]{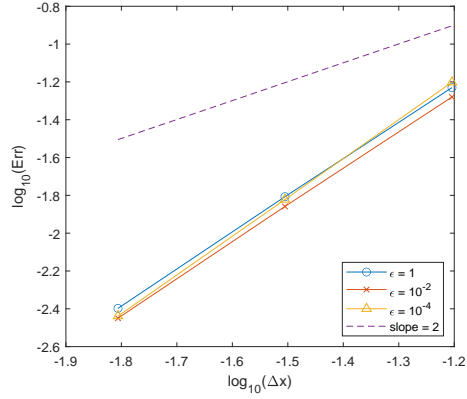}\\
		\caption{The test of convergence order with initial data \eqref{data0}. This figure shows the $l^2$ norm of fluid velocity $u$ with different $\varepsilon$.}\label{fig:ThreeEx0u}
	\end{figure}

	The convergence order in $l^{1}$ norm for the particle distribution $f_i (i=1, 2)$ and in $l^{2}$ norm for the fluid velocity $\mathbf{u}$ is reported in Figures \ref{fig:ThreeEx0f}-\ref{fig:ThreeEx0u}.
	This shows that the scheme is second order in space (hence in time) uniformly in $\varepsilon$ for both the particle distributions $f_i$ and the fluid velocity $\mathbf{u}$ as expected.

	\begin{figure}[htbp]
		\centering
		\includegraphics[width=15cm]{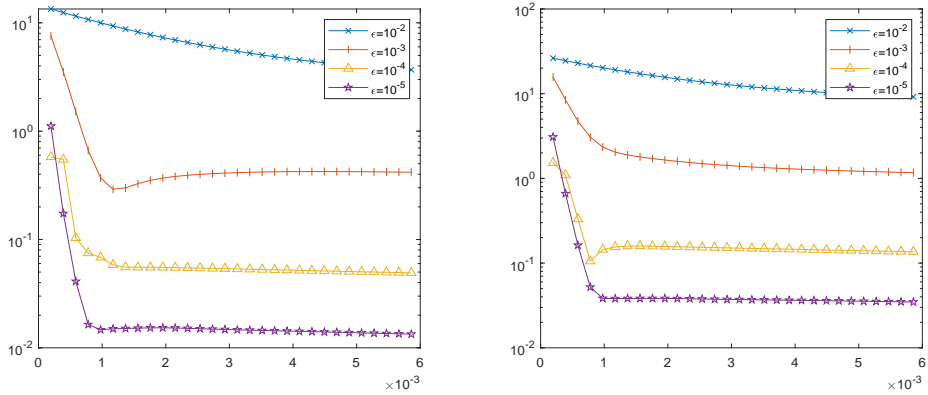}\\
		\caption{The time evolution $||f_i - n_i M_{u,i}||_1$ with different $\varepsilon$, starting with the initial data \eqref{Ex1Three}. The left is for the first particle $(i = 1)$ and the right is for the second particle  $(i = 2)$.}\label{fig:ex1err}
	\end{figure}

	\subsection{AP property}
	
	Now we check the AP property. We take the volcano like initial data: 
	\begin{equation}\label{Ex1Three}
		\begin{aligned}
			&n_1(0, \mathbf{x})=n_2(0, \mathbf{x})\\
			&\quad\quad\quad\  =\left(0.5+100\left((x_1-0.5)^{2}+(x_2-0.5)^{2}\right)\right) \exp \left(-40(x_1-0.5)^{2}-40(x_2-0.5)^{2}\right), \\
			&\mathbf{u}_{p,1}(0, \mathbf{x})=\mathbf{u}_{p,2}(0, \mathbf{x})=\left(\begin{array}{c}
				-\sin (2 \pi(x_2-0.5)) \\
				\sin (2 \pi(x_1-0.5))
			\end{array}\right) \exp \left(-20(x_1-0.5)^{2}-20(x_2-0.5)^{2}\right), \\
			&\mathbf{u}(0, \mathbf{x})=0.
		\end{aligned}
	\end{equation}
	
	The time evolution of $\ell^1$ distances $\left\|f_i-n_i M_{u,i}\right\|_1$ where $M_{u,i}$ is a Maxwellian centered at the fluid velocity $u$ is shown in Fig. \ref{fig:ex1err}. 
	The result gives a direct evidence of the AP property we proposed.

	
	\begin{figure}[htbp]
		\centering
		\includegraphics[width=13cm]{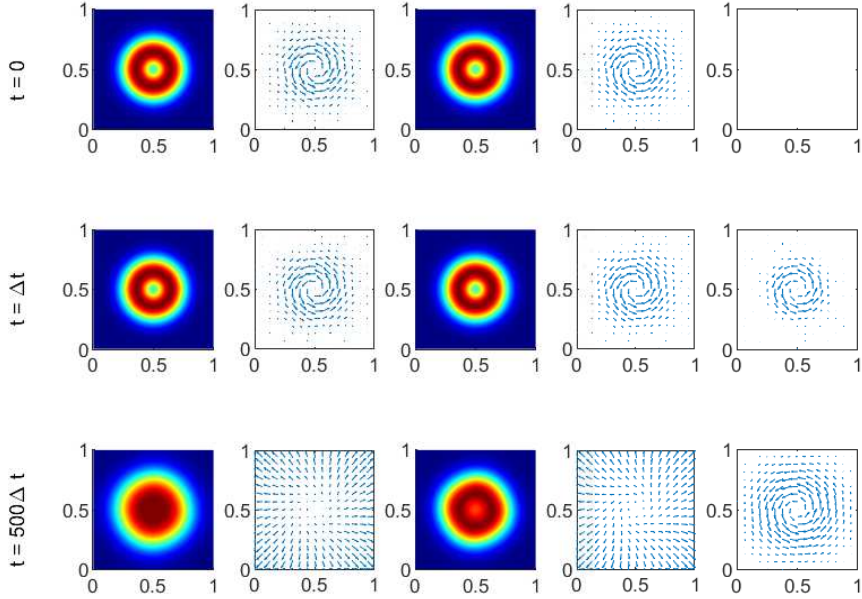}\\
		\caption{Numerical simulation with initial data \eqref{Ex1Three} for $\varepsilon=1$. The figure shows the particle density $n_{1}, n_{2}$ (first and third column), streamlines of particle velocity $\mathbf{u}_{p,1}, \mathbf{u}_{p,2}$ (second and fourth column) and fluid velocity $\mathbf{u}$ (fifth column) at $t=0$ (upper row), $t=\Delta t$ (middle row) and $t=500\Delta t$ (lower row).}\label{fig:ex1e1}
	\end{figure}
	
	\begin{figure}[htbp]
		\centering
		\includegraphics[width=13cm]{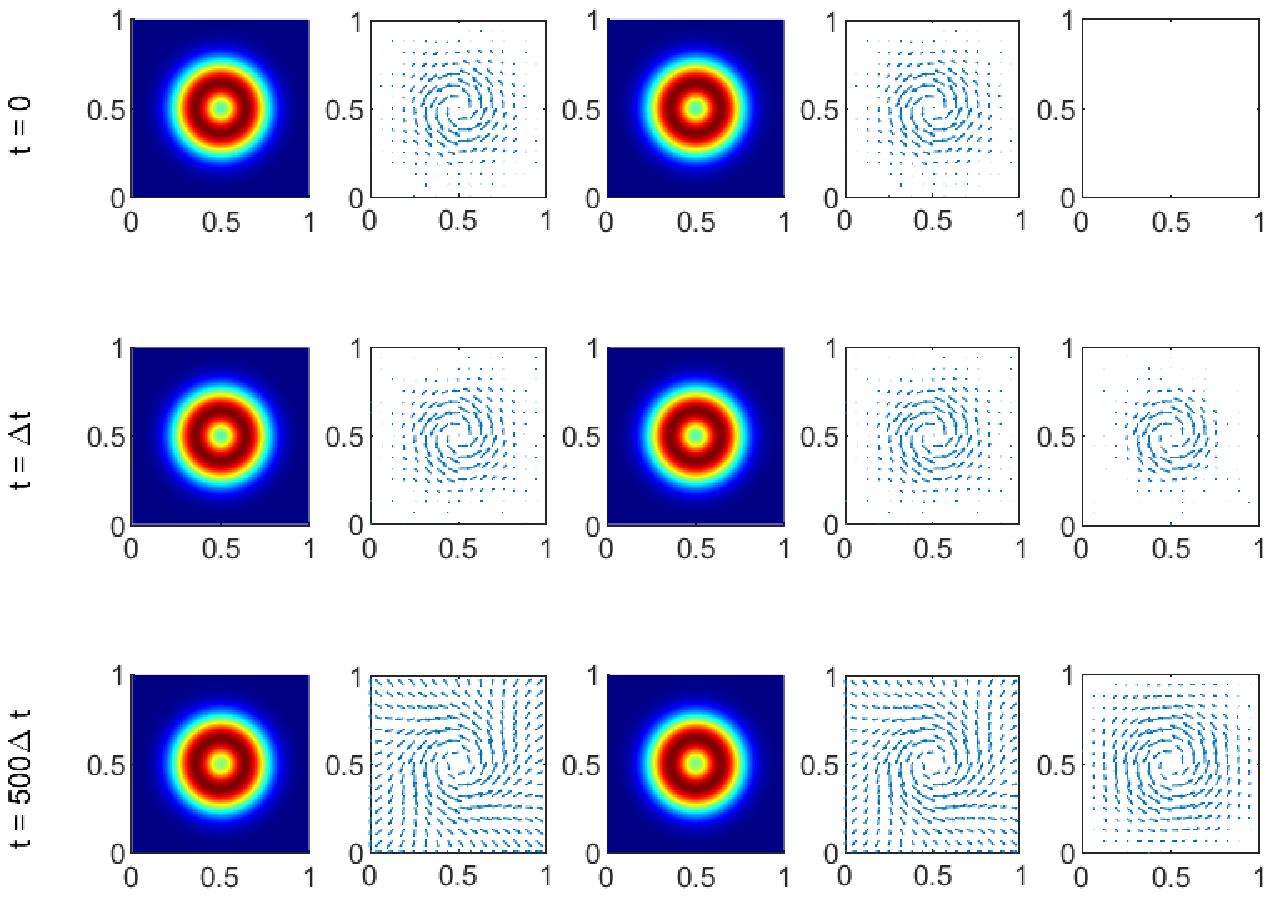}\\
		\caption{Numerical simulation with initial data \eqref{Ex1Three} for $\varepsilon=10^{-3}$. The figure shows the particle density $n_{1}, n_{2}$ (first and third column), streamlines of particle velocity $\mathbf{u}_{p,1}, \mathbf{u}_{p,2}$ (second and fourth column) and fluid velocity $\mathbf{u}$ (fifth column) at $t=0$ (upper row), $t=\Delta t$ (middle row) and $t=500\Delta t$ (lower row).}\label{fig:ex1e2}
	\end{figure}
	
	\begin{figure}[htbp]
		\centering
		\includegraphics[width=13cm]{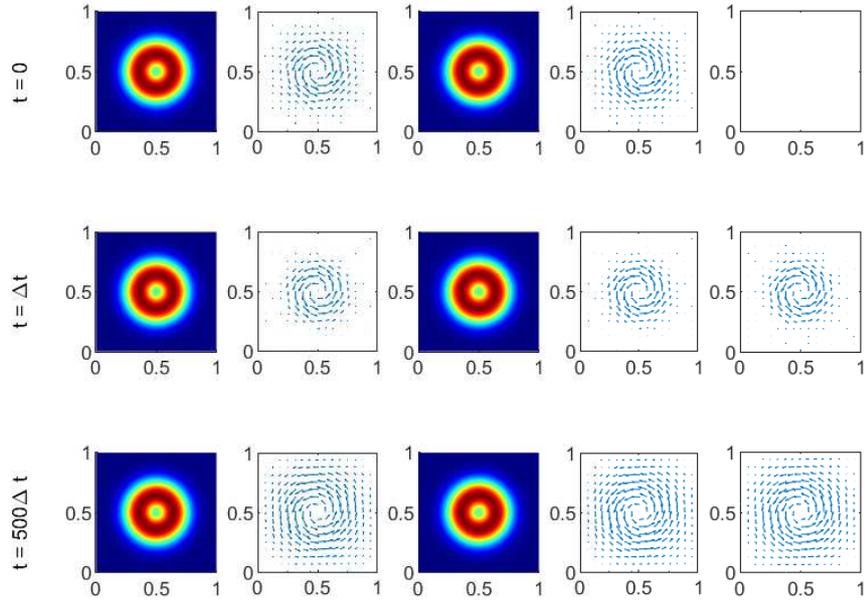}\\
		\caption{Numerical simulation with initial data \eqref{Ex1Three} for $\varepsilon=10^{-5}$. The figure shows the particle density $n_{1}, n_{2}$ (first and third column), streamlines of particle velocity $\mathbf{u}_{p,1}, \mathbf{u}_{p,2}$ (second and fourth column) and fluid velocity $\mathbf{u}$ (fifth column) at $t=0$ (upper row), $t=\Delta t$ (middle row) and $t=500\Delta t$ (lower row).}\label{fig:ex1e3}
	\end{figure}

	For the initial data \eqref{Ex1Three}, $\mathbf{u}_{p,i} \neq \mathbf{u}$ $(i=1, 2)$ and thus the equilibrium is not assumed. Figure \ref{fig:ex1e1} gives the time evolution of the system \eqref{Ex1Three} with $\varepsilon=1$ at $t^0$ (the initial time), $t^1$ (after one time step) and $t^{500}$ (the end time).
	The left two columns show the particle density $n_1$ and streamlines of particle velocity $\mathbf{u}_{p,1}$ of the first particle $(i=1)$. The third and fourth columns show the particle density $n_2$ and streamlines of particle velocity $\mathbf{u}_{p,2}$ of the second particle $(i=2)$. The right column give the fluid velocity $\mathbf{u}$. In Figure \ref{fig:ex1e1}, the particles expand to the whole square domain and are not significantly affected by the circulating fluid.
	The streamlines of particles and fluid are quite different because the drag force bewteen the fluid phase and particle phases is not significant. The behavior turns out to be significantly different as $\varepsilon$ decreases.
	
	Figures \ref{fig:ex1e2} gives the time evolution of this system \eqref{Ex1Three} with $\varepsilon=10^{-3}$. In the case with smaller $\varepsilon$, the drag force between different phases is much stronger that the particles also circulate in the square domain. Besides, the expansion in particle density is decelerated by the fluid. 
	
	Figures  \ref{fig:ex1e3} gives the time evolution of this system \eqref{Ex1Three} with  $\varepsilon=10^{-5}$. In this case, the particles stop expanding immediately due to the strong drag force.  The particles keep the volcano shape well in this period of time.
	

	\subsection{Some applications }
	\label{sec:num.3}
	
	In this section we apply our schemes to several different problems.
	In the following simulation we will take $Re=1000$. Larger Reynolds number, which requires smaller mesh size $\Delta x$ for the sake of accuracy, is beyond the scope of this work.
	In Section \ref{sec:num.1} the external force (the gravity) is considered and in Section \ref{sec:num.2} we perform a simulation with $\varepsilon$ varying in space.

	\subsubsection{Simulationof gravity driven flow}
	\label{sec:num.1}
	
	Now we consider the dam like initial data,
	\begin{equation}\label{Ex2Three}
		\begin{aligned}
			&n_1(0, \mathbf{x}) = n_2(0, \mathbf{x}) = 10^{-10}+1_{0 \leqslant x \leqslant 0.5}, \\
			&\mathbf{u}_{p,1}(0, \mathbf{x})= \mathbf{u}_{p,2}(0, \mathbf{x}) = 0, \\
			&\mathbf{u}(0, \mathbf{x})=0.
		\end{aligned}
	\end{equation}
	{\tiny }
	\begin{figure}[htbp]
		\centering
		\includegraphics[width=13cm]{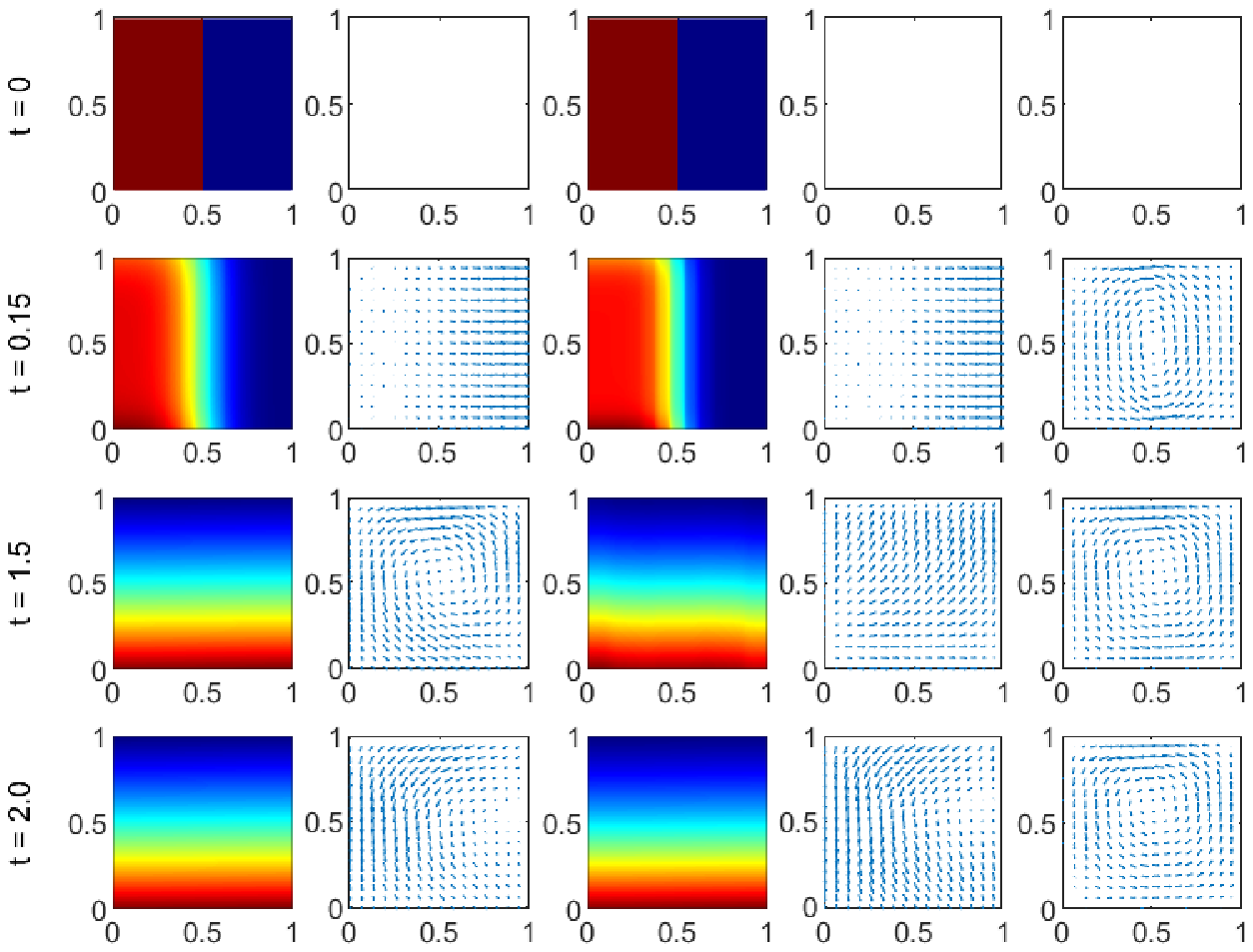}\\
		\caption{Numerical simulation with initial data \eqref{Ex2Three} for $\varepsilon=1$. The figure shows the particle density $n_{1}, n_{2}$ (first and third column), streamlines of particle velocity $\mathbf{u}_{p,1}, \mathbf{u}_{p,2}$ (second and fourth column) and fluid velocity $\mathbf{u}$ (fifth column) at $t=0, 0.15, 1.5, 2.0$.}\label{fig:ex2e1}
	\end{figure}
	
	\begin{figure}[htbp]
		\centering
		\includegraphics[width=14cm]{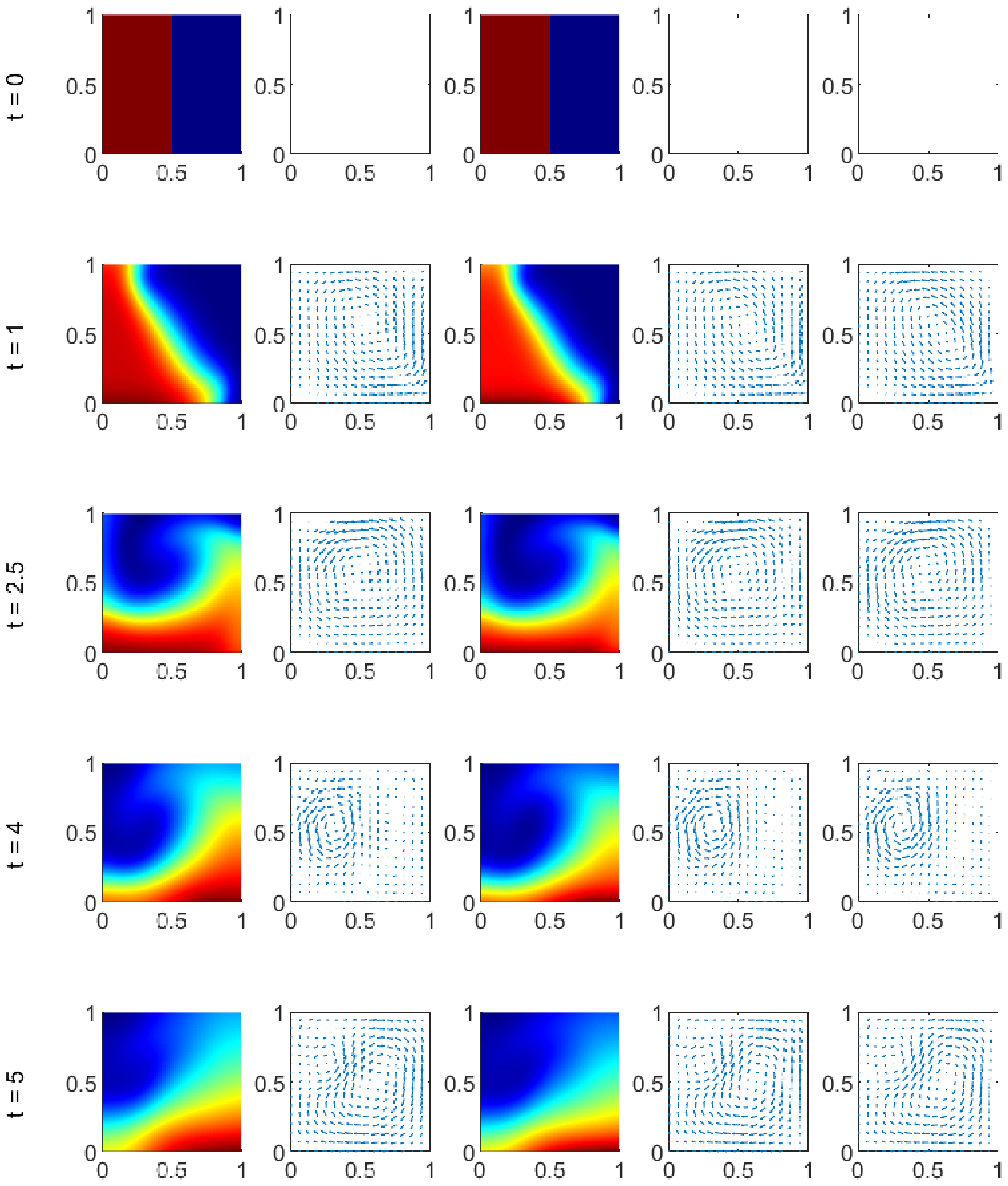}\\
		\caption{Numerical simulation with initial data \eqref{Ex2Three} for $\varepsilon=10^{-2}$. The figure shows the particle density $n_{1}, n_{2}$ (first and third column), streamlines of particle velocity $\mathbf{u}_{p,1}, \mathbf{u}_{p,2}$ (second and fourth column) and fluid velocity $\mathbf{u}$ (fifth column) at $t=0, 1, 2.5, 4, 5$.}\label{fig:ex2e2}
	\end{figure}
	
	\begin{figure}[htbp]
		\centering
		\includegraphics[width=14cm]{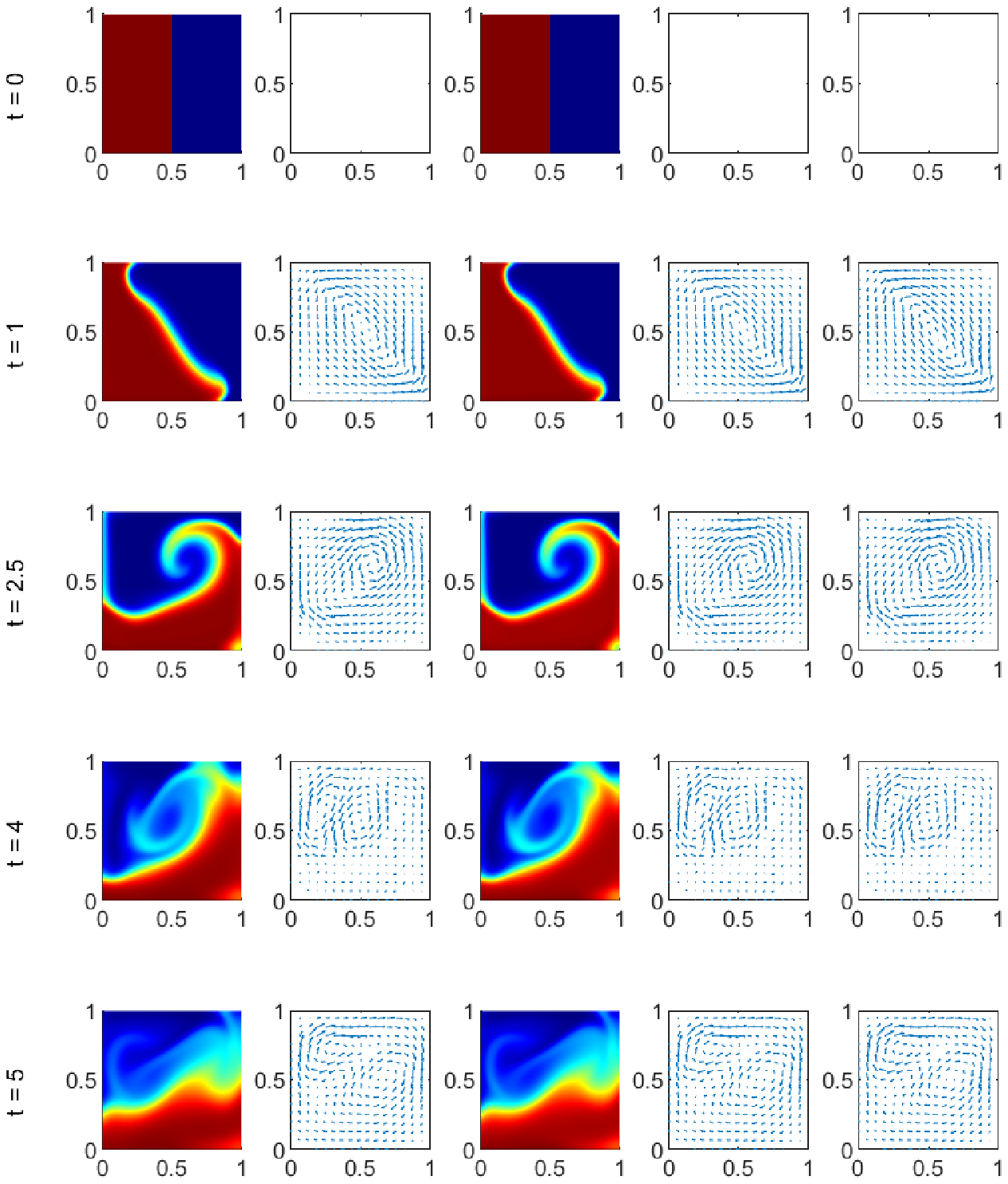}\\
		\caption{Numerical simulation with initial data \eqref{Ex2Three} for $\varepsilon=10^{-3}$. The figure shows the particle density $n_{1}, n_{2}$ (first and third column), streamlines of particle velocity $\mathbf{u}_{p,1}, \mathbf{u}_{p,2}$ (second and fourth column) and fluid velocity $\mathbf{u}$ (fifth column) at $t=0,1, 2.5, 4, 5$.}\label{fig:ex2e3}
	\end{figure}

	In this case the movement of particles and fluid are initiated by the gravity.  As the simulation starts, the particles fall down and cause the circulation of fluid.
	
	Fig. \ref{fig:ex2e1} shows the time evolution of the density and streamlines of the velocity for the small particle (left two columns) with $\varepsilon=1$, as well as the time evolution of  the density and streamlines of the velocity for the large particle (third and fourth column) and fluid velocity (right column). In this case the drag force between particles and fluid is not significant. The particles just
	fall down and cover the whole bottom. Clearly, the streamlines of particles is quite different from that of the fluid. 
	The first particle with smaller size evolves faster than the second particle with larger size. In other words, the light species gets close to the Maxwellian faster than the heavy one, which is consistent with the observation verified in \cite{JinLi2013} for multi-species Boltzmann equations. 
	
	Fig. \ref{fig:ex2e2} shows the time evolution of particle density with $\varepsilon=10^{-2}$. Now the drag force between particles and fluid is stronger. 
	The streamlines of particles and the fluid are similar all the time. As time evolves, the particles fall down and drive the fluid to circulate counter-clockwisely. Then the particles follow this circulation. Finally the particles settle at the bottom uniformly due to the loss of energy.
	
	Fig. \ref{fig:ex2e3} shows the time evolution of particle density with $\varepsilon=10^{-3}$. 
	In this case the drag force between particles and fluid is so strong that the streamlines of the particles and fluid are quite similar to each other.

	\subsubsection{Simulation of injecting problem}
	\label{sec:num.2}
	
\begin{figure}[htbp]
		\centering
		\includegraphics[width=8cm]{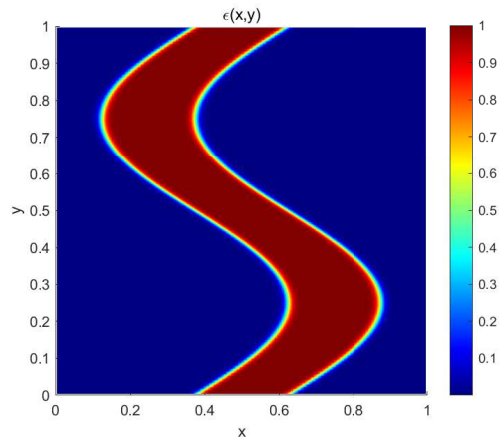}\\
		\caption{The $\mathbf{x}$-dependent function $\varepsilon(\mathrm{x})$ given by \eqref{Ex30}.
		}\label{fig:Ex4eps}
	\end{figure}

	One of the advantages of AP schemes is that they can capture the solution behaviors automatically as $\varepsilon$ varies in space. Finally, let us  consider a mixing regime problem, with an $\mathbf{x}$-dependent $\varepsilon(\mathbf{x})$
	\begin{equation}\label{Ex30}
		\varepsilon(x, y)=\varepsilon_0+\frac{1}{2}\left(\tanh \left(10-80\left(x-\frac{1}{2}-\frac{1}{4} \sin (2 \pi y)\right)\right)+\tanh \left(10+80\left(x-\frac{1}{2}-\frac{1}{4} \sin (2 \pi y)\right)\right)\right).
	\end{equation}
	Here $\varepsilon_0 \ll 1$ is a constant. $\varepsilon(\mathbf{x})$ varies from $\varepsilon_0$ to $O(1)$ smoothly, as shown in Figure \ref{fig:Ex4eps} with $\varepsilon_0=10^{-3}$.

	\begin{figure}[htbp]
		\centering
		\includegraphics[width=14cm]{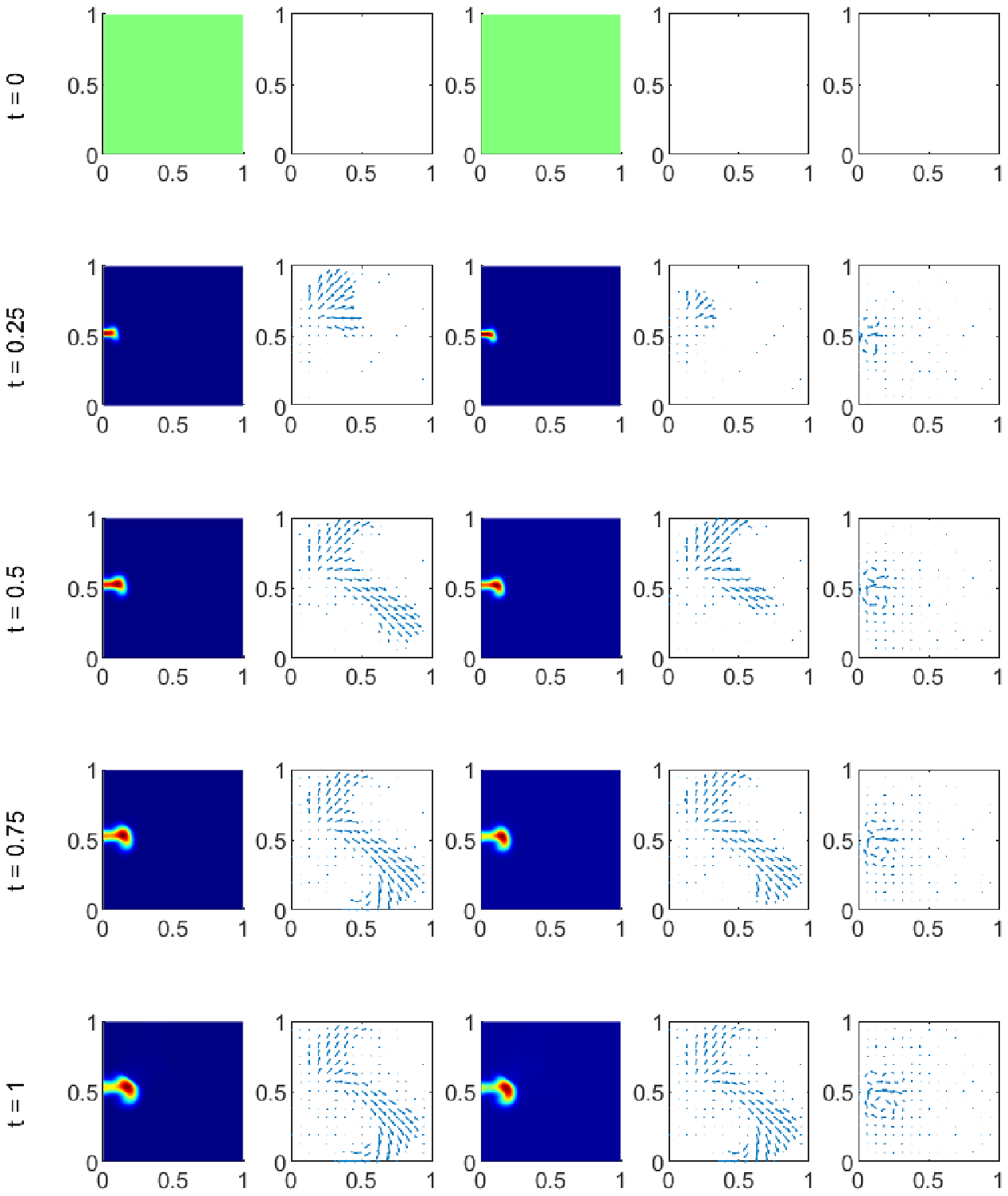}\\
		\caption{Numerical simulation with initial data \eqref{Ex30}-\eqref{Ex3Three2}. The figure shows the particle density $n_{1}, n_{2}$ (first and third column), streamlines of particle velocity $\mathbf{u}_{p,1}, \mathbf{u}_{p,2}$ (second and fourth column) and fluid velocity $\mathbf{u}$ (fifth column) at $t=0, 0.5, 2.5, 5$.}\label{fig:Ex4}
	\end{figure}

	Consider the situation when the particles are injected into the square domain. 
	More specifically, the initial conditions are given as follows.
	\begin{equation}\label{Ex3Three}
		n_1(0, \mathbf{x})=n_2(0, \mathbf{x})=10^{-10}, \quad \mathbf{u}_{p,1}(0, \mathbf{x})=\mathbf{u}_{p,2}(0, \mathbf{x})=\mathbf{u}(0, \mathbf{x})=0 .
	\end{equation}
	The injecting particle flow is described by the boundary condition on $f_i (i=1,2)$,
	\begin{equation}\label{Ex3Three2}
		\begin{aligned}
			&f_1(t, \mathbf{x}, \mathbf{v})=1_{2 \leqslant v_{1} \leqslant 3}, \quad \text { if } \mathbf{x} \in \Gamma = \{(0, y) \mid 0.475 \leqslant y \leqslant 0.575\}, \\
			&f_2(t, \mathbf{x}, \mathbf{v})=1_{2 \leqslant v_{1} \leqslant 3}, \quad \text { if } \mathbf{x} \in \Gamma = \{(0, y) \mid 0.45 \leqslant y \leqslant 0.55\},
		\end{aligned}	
	\end{equation}
	where $v_{1}$ is the first component of $\mathbf{v}$. The entrance of flow $\Gamma$ locates at the center of the left boundary.


	Figure \ref{fig:Ex4} shows several snapshots of the time evolution of the particle density $n_1, n_2$, streamlines of prticle velocity $u_{p,1}, u_{p,2}$ and the fluid velocity $u$. The behaviors of both phases are clearly influenced by the spatially variable $\varepsilon(\mathbf{x})$. 
	The first particle with smaller size evolves a little faster than the second particle with larger size. 
	The difference of the velocities of two phases $\left|\mathbf{u}_p-\mathbf{u}\right|$ shows an S-shape profile which is consistent to $\varepsilon(\mathbf{x})$ in Figure \ref{fig:Ex4eps}. This suggests that the fluid limit of this two-phase system is achieved automatically in the strong interaction regime where $\varepsilon \ll 1$. While in the weak interaction regimes where $\varepsilon=O(1)$ the two phases behave quite differently. 

	\section{Conclusion}
	\label{sec:conclusion}
	We have introduced a new numerical scheme for simulating a system coupling the incompressible Navier–Stokes equation to multi-component Vlasov–Fokker–Planck equations with distinct particle sizes. In particular, the scheme is Asymptotic Preserving, which drives the system towards a set of hydrodynamic equations in the regime of $\varepsilon \ll 1$. The method is based on the implicit treatment of the potentially stiff terms and relies on the possibility of updating the macroscopic unknowns by solving simple linear systems. 
	It is interesting that in the multi-phase kinetic-fluid system, the smaller particle with light mass evolves a little faster than the bigger particle with heavy mass, and thus the smaller particle may get close to the equilibrium faster than the bigger one.
	The scheme can easily incorporate relevant generalizations of the basic model like considering variable fluid density or temperature-dependent viscosities.
	Furthermore, the scheme can be used for the kinetic-fluid system with multi-size particles with uncertainties.

	\section*{Acknowledgement}
	S. Jin was partially supported by National Key R\&D Program of China (no. 2020YFA0712000) and National Natural Science Foundation of China (no. 20Z103020029). Y. Lin was partially supported by National Natural Science Foundation of China (no. 12201404) and Postdoctoral Research Foundation of China (no. 2021M702142 and no. 2021TQ0203). Y. Lin thanks Dr.  Ruiwen Shu for his  help during the preparation of the paper.

	\section*{Appendix}
	 \textbf{Scaling issues.} Let us detail precisely the scaling issues and the physical meaning of the system we are interested in. 
		To this end, let us go back to the equations written with dimensional quantities. 
	
	We adopt a multi-component Vlasov-Fokker-Planck-incompressible Navier-Stokes system for a discrete modeling of $N$-size variable for  particles. Let $i \in \{1, 2,\ldots N\}$. We refer to " a particle size $i$ " as an assembly of $i$ monomers. 
	As usual, the fluid is described by its density $n(t, x)$ and velocity $u(t, x)$, which obey the mass conservation and momentum balance relations.
	Denoting by $a>0$ the  radius of a monomer and $\rho_P$ its mass density, the radius and mass of an $i$-mer are defined by $$r_i=a i^{1 / 3}, \quad m_i=\frac{4}{3} \pi a^3 i \rho_P = m_1 i \text{~~with~} m_1 = \frac{4}{3} \pi a^3  \rho_P$$
	respectively.
	The particles are described by the quantity $f_i(t,x,v)$ such that $f_i(t,x,v) \mathrm{~d}v \mathrm{~d}x$ is the probability of finding particles with size $i$ in the domain of the phase space centered at $(x,v)$ with volume $\mathrm{~d}v \mathrm{~d}x$.
	Particles are subject to the Brownian motion, which induces a diffusion term whose coefficient is given by Einstein formula \cite{Einstein19}
	$$
	\frac{9 \mu}{2 \rho_P r_i^2} \frac{k \theta}{m_i} 
	$$
	where $k$ is the Boltzmann constant, $\theta>0$ is temperature of the fluid and $\mu$ is dynamic viscosity of the fluid. 
	The drag force exerted by the fluid on the particles is given by the 
	Stokes law, which is proportional to the relative velocity
	$$6 \pi \mu r_i(v-u)=6 \pi \mu a i^{1 / 3}(v-u).$$
	The Sokes settling time is defined by
	$$
	\tau_i=\frac{m_i}{6 \pi \mu r_i}=\dfrac{2 \rho_P r_i^2}{9 \mu}=i^{2 / 3} \tau_1 \quad \text { with } \tau_1=\dfrac{2 \rho_P a^2}{9 \mu},
	$$
	which is typical of the effect of the drag force on the $i$-particle.
	Denoting by $\rho_F>0$ a typical mass density for the fluid and introducing the $i$ th Fokker-Planck (FP) operator $\mathcal{L}_u f_i$ by
	$$\mathcal{L}_u f_i=\nabla_v \cdot\left((v-u) f_i+\dfrac{k \theta}{m_i} \nabla_v f_i\right),$$
	the particles-fluid mixture with distinct particle size is then described by the followign system of PDEs:
	\begin{equation}\label{ModelPDEs}
		\left\{\begin{array}{l}\partial_t f_i+v \cdot \nabla_x f_i-\nabla_x \Phi \cdot \nabla_v f_i=\dfrac{9 \mu}{2 \rho_P r_i^2} \mathcal{L}_u f_i, \\ \rho_F\partial_t(u)+\rho_F\operatorname{Div}_x(u \otimes u)+\rho_F\alpha \nabla_x \Phi+ \nabla_x p-\mu \Delta_x u=\displaystyle 6 \pi \mu \sum_{i=1}^{N} \int_{\mathbb{R}^3}(v-u) f_i r_i \mathrm{~d} v, \\ \nabla_x \cdot u=0.\end{array}\right.
	\end{equation}


	Following \cite{Carrilo2006,Goudon2013fluid}, we are going to write system \eqref{ModelPDEs} in dimensionless form. 
	To this end, we introduce time and length scales, denoted by $T$ and $L$, respectively, which define the velocity unit $U=L / T$.
	Set the thermal velocity
	$V=\sqrt{\frac{k \bar{\theta}}{m_1}}$ with $m_1=\frac{4}{3} \pi a^3 \rho_P$ and $\bar{\theta}>0$ a reference temperature. Set $\mathcal{P}$ a suitable pressure unit.
	The dimensionless variables and unknowns can be defined as follows:
	\begin{itemize}
		\item $t^{\prime}=t/T, x^{\prime}=x/L,  v^{\prime}=v/V$,
		\item $n^{\prime}\left(t^{\prime}, x^{\prime}\right)=n\left(T t^{\prime}, L x^{\prime}\right),  u^{\prime}\left(t^{\prime}, x^{\prime}\right)=u\left(T t^{\prime}, L x^{\prime}\right)/U$,
		\item $ p^{\prime}\left(t^{\prime}, x^{\prime}\right)=p\left(T t^{\prime}, L x^{\prime}\right)/\mathcal{P}, f_i^{\prime}\left(t^{\prime}, x^{\prime}, v^{\prime}\right)=\frac{4}{3} \pi a^3 V^3 f_i\left(T t^{\prime}, L x^{\prime}, V v^{\prime}\right)$.
	\end{itemize}  
	Since $ \mathrm{~d} v^{\prime}=\mathrm{d} v/V^3$, for any given function $\varphi$, one has
	$$
	\int_{\mathbb{R}^3} \varphi(v) f_i(t, x, v) \mathrm{d} v=\frac{1}{\frac{4}{3} \pi a^3} \int_{\mathbb{R}^3} \varphi\left(V v^{\prime}\right) f_i^{\prime}\left(t^{\prime}, x^{\prime}, v^{\prime}\right) \mathrm{d} v^{\prime}.
	$$
	If the temperature is not assumed constant, similarly set $ \theta^{\prime}\left(t^{\prime}, x^{\prime}\right)=\theta\left(T t^{\prime}, L x^{\prime}\right)/\bar{\theta}$. 
	For the external potential, set $ \Phi^{\prime}\left(x^{\prime}\right) = \Phi\left(L x^{\prime}\right)\frac{\tau_1}{\vartheta_s L}$, where $\vartheta_s$ has the dimension of velocity (for gravity-driven flows, it is the Stokes settling velocity). One arrives at 
	$$
	\begin{aligned}
		\frac{1}{T} \partial_{t^{\prime}}\left(f_i^{\prime}\right) & +\frac{V}{L} v^{\prime} \cdot \nabla_{x^{\prime}}\left(f_i^{\prime}\right)-\frac{\vartheta_s}{\tau_1 V} \nabla_{x^{\prime}} \Phi^{\prime} \cdot \nabla_{v^{\prime}}\left(f_i\right) 
		=\frac{1}{\tau_i} \nabla_{v^{\prime}} \cdot\left(\left(v^{\prime}-\frac{U}{V} u^{\prime}\right) f_i^{\prime}+\frac{k \bar{\theta}}{m_i V^2} \theta^{\prime} \nabla_{v^{\prime}} f_i^{\prime}\right).
	\end{aligned}
	$$ 
	Therefore, we realize that the system is driven by the following set of dimensionless paratmeters:
	$$
	\beta=\frac{T}{L} V=\frac{V}{U}, \quad \frac{1}{\varepsilon}=\frac{T}{\tau_1}, \quad \eta=\frac{\vartheta_s T}{V \tau_1}, \quad \chi=\frac{\mathcal{P} T}{\rho_F L U}=\frac{\mathcal{P}}{\rho_F U^2},
	$$
	together with the density ratio $\rho_P / \rho_F$. Finally, by dropping the prime marks, one obtains the dimensionles form of system \eqref{ModelPDEs}:
	\begin{equation}\label{ModelScaled}
		\begin{aligned}
			& \left\{\begin{array}{l}
				\partial_t f_i+\beta v \cdot \nabla_x f_i-\eta \nabla_x \Phi \cdot \nabla_v f_i=\dfrac{1}{\varepsilon} \dfrac{1}{i^{2 / 3}} \nabla_v \cdot\left(\left(v-\dfrac{1}{\beta} u\right) f_i+\dfrac{\bar{\theta}}{i} \nabla_v f_i\right), \\
				\partial_t(u)+\operatorname{Div}_x(u \otimes u)+\alpha \beta \eta \nabla_x \Phi+\chi \nabla_x p=\dfrac{\rho_P}{\varepsilon \rho_F} \displaystyle\sum_{i=1}^{N} \displaystyle\int_{\mathbb{R}^3}(\beta v-u) f_i i^{1 / 3} \mathrm{~d} v+\mu \Delta_x u, \\
				\nabla_x \cdot u=0.
			\end{array}\right. 
		\end{aligned}
	\end{equation}
	Here $\mu$ stands for the rescaled and dimensionless version of the fluid viscosity.

	\bibliography{mybibfile}

\begin{thebibliography}{37}
\expandafter\ifx\csname natexlab\endcsname\relax\def\natexlab#1{#1}\fi
\providecommand{\url}[1]{\texttt{#1}}
\providecommand{\href}[2]{#2}
\providecommand{\path}[1]{#1}
\providecommand{\DOIprefix}{doi:}
\providecommand{\ArXivprefix}{arXiv:}
\providecommand{\URLprefix}{URL: }
\providecommand{\Pubmedprefix}{pmid:}
\providecommand{\doi}[1]{\href{http://dx.doi.org/#1}{\path{#1}}}
\providecommand{\Pubmed}[1]{\href{pmid:#1}{\path{#1}}}
\providecommand{\bibinfo}[2]{#2}
\ifx\xfnm\relax \def\xfnm[#1]{\unskip,\space#1}\fi
\bibitem[{Andrews and O'Rourke(1996)}]{Andrews1996}
\bibinfo{author}{Andrews, M.J.}, \bibinfo{author}{O'Rourke, P.J.},
  \bibinfo{year}{1996}.
\newblock \bibinfo{title}{The multiphase particle-in-cell (mp-pic) method for
  dense particulate flows}.
\newblock \bibinfo{journal}{International Journal of Multiphase Flow}
  \bibinfo{volume}{22}, \bibinfo{pages}{379--402}.
\bibitem[{Aregba-Driollet and Milisik(2004)}]{Aregba2004}
\bibinfo{author}{Aregba-Driollet, D.}, \bibinfo{author}{Milisik, V.},
  \bibinfo{year}{2004}.
\newblock \bibinfo{title}{Kinetic approximation of a boundary value problem for
  conservation laws}.
\newblock \bibinfo{journal}{Numerische Mathematik} \bibinfo{volume}{97},
  \bibinfo{pages}{595–633}.
\bibitem[{Baranger et~al.(2005)Baranger, Boudin, Jabin and
  Mancini}]{Baranger2005}
\bibinfo{author}{Baranger, C.}, \bibinfo{author}{Boudin, L.},
  \bibinfo{author}{Jabin, P.E.}, \bibinfo{author}{Mancini, S.},
  \bibinfo{year}{2005}.
\newblock \bibinfo{title}{A modeling of biospray for the upper airways, in
  cemracs 2004—mathematics and applications to biology and medicine}.
\newblock \bibinfo{journal}{ESAIM: Proc. 14, EDP Sci., Les Ulis, France} ,
  \bibinfo{pages}{41--47}.
\bibitem[{Baranger and Desvillettes(2006)}]{Baranger2006}
\bibinfo{author}{Baranger, C.}, \bibinfo{author}{Desvillettes, L.},
  \bibinfo{year}{2006}.
\newblock \bibinfo{title}{Coupling euler and vlasov equations in the context of
  sprays: the local-in-time, classical solutions}.
\newblock \bibinfo{journal}{Journal of Hyperbolic Differential Equations}
  \bibinfo{volume}{3}, \bibinfo{pages}{1--26}.
\bibitem[{Boudin et~al.(2009a)Boudin, Boutin, Fornet, Goudon, Lafitte,
  Lagouti{\`e}re and Merlet}]{Boudin2009fluid}
\bibinfo{author}{Boudin, L.}, \bibinfo{author}{Boutin, B.},
  \bibinfo{author}{Fornet, B.}, \bibinfo{author}{Goudon, T.},
  \bibinfo{author}{Lafitte, P.}, \bibinfo{author}{Lagouti{\`e}re, F.},
  \bibinfo{author}{Merlet, B.}, \bibinfo{year}{2009}a.
\newblock \bibinfo{title}{Fluid-particles flows: A thin spray model with energy
  exchanges}, in: \bibinfo{booktitle}{ESAIM: Proceedings},
  \bibinfo{organization}{EDP Sciences}. pp. \bibinfo{pages}{195--210}.
\bibitem[{Boudin et~al.(2009b)Boudin, Desvillettes, Grandmont and
  Moussa}]{BoudinDesvillettes2009}
\bibinfo{author}{Boudin, L.}, \bibinfo{author}{Desvillettes, L.},
  \bibinfo{author}{Grandmont, C.}, \bibinfo{author}{Moussa, A.},
  \bibinfo{year}{2009}b.
\newblock \bibinfo{title}{Global existence of solutions for the coupled vlasov
  and navier-stokes equations}.
\newblock \bibinfo{journal}{Differential and Integral Equations}
  \bibinfo{volume}{22}, \bibinfo{pages}{1247--1271}.
\bibitem[{Carrillo et~al.(2011)Carrillo, Duan and Moussa}]{CarrilloDuan2011}
\bibinfo{author}{Carrillo, J.A.}, \bibinfo{author}{Duan, R.},
  \bibinfo{author}{Moussa, A.}, \bibinfo{year}{2011}.
\newblock \bibinfo{title}{Global classical solutions close to equilibrium to
  the vlasov-fokker-planck-euler system}.
\newblock \bibinfo{journal}{Kinetic \& Related Models} \bibinfo{volume}{4},
  \bibinfo{pages}{227--258}.
\bibitem[{Carrillo et~al.(2008)Carrillo, Goudon and Lafitte}]{Carrillo2008}
\bibinfo{author}{Carrillo, J.A.}, \bibinfo{author}{Goudon, T.},
  \bibinfo{author}{Lafitte, P.}, \bibinfo{year}{2008}.
\newblock \bibinfo{title}{Simulation of fluid and particles flows: Asymptotic
  preserving schemes for bubbling and flowing regimes}.
\newblock \bibinfo{journal}{Journal of Computational Physics}
  \bibinfo{volume}{227}, \bibinfo{pages}{7929--7951}.
\bibitem[{Carrilo and Goudon(2006)}]{Carrilo2006}
\bibinfo{author}{Carrilo, J.}, \bibinfo{author}{Goudon, T.},
  \bibinfo{year}{2006}.
\newblock \bibinfo{title}{Stability and asymptotic analysis of a
  fluid-particles interaction model}.
\newblock \bibinfo{journal}{Communications in Partial Differential Equations}
  \bibinfo{volume}{31}, \bibinfo{pages}{1349--1379}.
\bibitem[{Chorin(1967)}]{Chorin1967}
\bibinfo{author}{Chorin, A.J.}, \bibinfo{year}{1967}.
\newblock \bibinfo{title}{The numerical solution of the navier-stokes equations
  for an incompressible fluid}.
\newblock \bibinfo{journal}{Bulletin of the American Mathematical Society}
  \bibinfo{volume}{73}, \bibinfo{pages}{928--931}.
\bibitem[{Chorin(1969)}]{Chorin1969}
\bibinfo{author}{Chorin, A.J.}, \bibinfo{year}{1969}.
\newblock \bibinfo{title}{On the convergence of discrete approximations to the
  navier-stokes equations}.
\newblock \bibinfo{journal}{Mathematics of computation} \bibinfo{volume}{23},
  \bibinfo{pages}{341--353}.
\bibitem[{Cober and Isaac(2006)}]{Cober2006}
\bibinfo{author}{Cober, S.}, \bibinfo{author}{Isaac, G.}, \bibinfo{year}{2006}.
\newblock \bibinfo{title}{Estimating Maximum Aircraft Icing Environments Using
  a Large Database of In-Situ Observations}.
\newblock 44th AIAA Aerospace Sciences Meeting and Exhibit,
  \bibinfo{publisher}{AIAA 2006-266}.
\bibitem[{Einstein(1906)}]{Einstein19}
\bibinfo{author}{Einstein, A.}, \bibinfo{year}{1906}.
\newblock \bibinfo{title}{Eine neue bestimmung der molek\"uldimensionen}.
\newblock \bibinfo{journal}{Ann. Physik} \bibinfo{volume}{19},
  \bibinfo{pages}{289--306}.
\bibitem[{Friedlander(1977)}]{Friedlander1977}
\bibinfo{author}{Friedlander, S.K.}, \bibinfo{year}{1977}.
\newblock \bibinfo{title}{Smoke, Dust and Haze: Fundamentals of Aerosol
  Behavior}.
\newblock \bibinfo{publisher}{Wiley-Interscience}, \bibinfo{address}{New York}.
\bibitem[{Goudon et~al.(2010)Goudon, He, Moussa and Zhang}]{GoudonHe2010}
\bibinfo{author}{Goudon, T.}, \bibinfo{author}{He, L.},
  \bibinfo{author}{Moussa, A.}, \bibinfo{author}{Zhang, P.},
  \bibinfo{year}{2010}.
\newblock \bibinfo{title}{The navier–stokes–vlasov–fokker–planck system
  near equilibrium}.
\newblock \bibinfo{journal}{SIAM Journal on Mathematical Analysis}
  \bibinfo{volume}{42}, \bibinfo{pages}{2177--2202}.
\bibitem[{Goudon et~al.(2004a)Goudon, Jabin and Vasseur}]{GoudonJabin2004a}
\bibinfo{author}{Goudon, T.}, \bibinfo{author}{Jabin, P.E.},
  \bibinfo{author}{Vasseur, A.}, \bibinfo{year}{2004}a.
\newblock \bibinfo{title}{{Hydrodynamic limit for the Vlasov-Navier-Stokes
  equations. I. Light particles regime}}.
\newblock \bibinfo{journal}{Indiana University Mathematics Journal}
  \bibinfo{volume}{53}, \bibinfo{pages}{1517--1536}.
\bibitem[{Goudon et~al.(2004b)Goudon, Jabin and Vasseur}]{GoudonJabin2004b}
\bibinfo{author}{Goudon, T.}, \bibinfo{author}{Jabin, P.E.},
  \bibinfo{author}{Vasseur, A.}, \bibinfo{year}{2004}b.
\newblock \bibinfo{title}{{Hydrodynamic limit for the Vlasov-Navier-Stokes
  equations. II. Fine particles regime}}.
\newblock \bibinfo{journal}{Indiana University Mathematics Journal}
  \bibinfo{volume}{53}, \bibinfo{pages}{1495--1515}.
\bibitem[{Goudon et~al.(2013a)Goudon, Jin, Liu and Yan}]{Goudon2013asymptotic}
\bibinfo{author}{Goudon, T.}, \bibinfo{author}{Jin, S.}, \bibinfo{author}{Liu,
  J.G.}, \bibinfo{author}{Yan, B.}, \bibinfo{year}{2013}a.
\newblock \bibinfo{title}{Asymptotic-preserving schemes for kinetic-fluid
  modeling of disperse two-phase flows}.
\newblock \bibinfo{journal}{Journal of Computational Physics}
  \bibinfo{volume}{246}, \bibinfo{pages}{145--164}.
\bibitem[{Goudon et~al.(2012)Goudon, Jin and Yan}]{Goudon2012}
\bibinfo{author}{Goudon, T.}, \bibinfo{author}{Jin, S.}, \bibinfo{author}{Yan,
  B.}, \bibinfo{year}{2012}.
\newblock \bibinfo{title}{Simulation of fluid-particles flows: heavy particles,
  flowing regime and asymptotic-preserving schemes}.
\newblock \bibinfo{journal}{Communications in Mathematical Sciences}
  \bibinfo{volume}{10}, \bibinfo{pages}{355--385}.
\bibitem[{Goudon et~al.(2013b)Goudon, Sy and Tine}]{Goudon2013fluid}
\bibinfo{author}{Goudon, T.}, \bibinfo{author}{Sy, M.}, \bibinfo{author}{Tine,
  L.M.}, \bibinfo{year}{2013}b.
\newblock \bibinfo{title}{A fluid-kinetic model for particulate flows with
  coagulation and breakup: stationary solutions, stability, and hydrodynamic
  regimes}.
\newblock \bibinfo{journal}{SIAM Journal on Applied Mathematics}
  \bibinfo{volume}{73}, \bibinfo{pages}{401--421}.
\bibitem[{Hamdache(1998)}]{Hamdache1998}
\bibinfo{author}{Hamdache, K.}, \bibinfo{year}{1998}.
\newblock \bibinfo{title}{Global existence and large time behaviour of
  solutions for the vlasov-stokes equations}.
\newblock \bibinfo{journal}{Japan Journal of Industrial and Applied
  Mathematics} \bibinfo{volume}{51}, \bibinfo{pages}{51--74}.
\bibitem[{Hauf and Schr\"oder(2006)}]{Hauf2006}
\bibinfo{author}{Hauf, T.}, \bibinfo{author}{Schr\"oder, F.},
  \bibinfo{year}{2006}.
\newblock \bibinfo{title}{Aircraft icing research flights in embedded
  convection}.
\newblock \bibinfo{journal}{Meteorology and Atmospheric Physics}
  \bibinfo{volume}{91}, \bibinfo{pages}{247--265}.
\bibitem[{Jin(1999)}]{Jin1999}
\bibinfo{author}{Jin, S.}, \bibinfo{year}{1999}.
\newblock \bibinfo{title}{{Efcient asymptotic-preserving (AP) schemes for some
  multiscale kinetic equations}}.
\newblock \bibinfo{journal}{SIAM Journal on Scientific Computing}
  \bibinfo{volume}{21}, \bibinfo{pages}{441--454}.
\bibitem[{Jin(2010)}]{Jin2010}
\bibinfo{author}{Jin, S.}, \bibinfo{year}{2010}.
\newblock \bibinfo{title}{Asymptotic preserving (ap) schemes for multiscale
  kinetic and hyperbolic equations: a review}.
\newblock \bibinfo{journal}{Lecture notes for summer school on methods and
  models of kinetic theory (M\&MKT), Porto Ercole (Grosseto, Italy)} ,
  \bibinfo{pages}{177--216}.
\bibitem[{Jin(2022)}]{Jin2022}
\bibinfo{author}{Jin, S.}, \bibinfo{year}{2022}.
\newblock \bibinfo{title}{Asymptotic-preserving schemes for multiscale physical
  problems}.
\newblock \bibinfo{journal}{Acta Numerica} \bibinfo{volume}{31},
  \bibinfo{pages}{415--489}.
\bibitem[{Jin and Li(2013)}]{JinLi2013}
\bibinfo{author}{Jin, S.}, \bibinfo{author}{Li, Q.}, \bibinfo{year}{2013}.
\newblock \bibinfo{title}{A bgk-penalization-based asymptotic-preserving scheme
  for the multispecies boltzmann equation}.
\newblock \bibinfo{journal}{Numerical Methods for Partial Differential
  Equations} \bibinfo{volume}{29}, \bibinfo{pages}{1056--1080}.
\bibitem[{Jin and Lin(2022)}]{JinLin2022}
\bibinfo{author}{Jin, S.}, \bibinfo{author}{Lin, Y.}, \bibinfo{year}{2022}.
\newblock \bibinfo{title}{Energy estimates and hypocoercivity analysis for a
  multi-phase navier-stokes-vlasov-fokker-planck system with uncertainty}
  \href{http://arxiv.org/abs/2204.10573}{{\tt arXiv:2204.10573}}.
\bibitem[{Jin and Yan(2011)}]{JinYan2011}
\bibinfo{author}{Jin, S.}, \bibinfo{author}{Yan, B.}, \bibinfo{year}{2011}.
\newblock \bibinfo{title}{A class of asymptotic-preserving schemes for the
  fokker--planck--landau equation}.
\newblock \bibinfo{journal}{Journal of Computational Physics}
  \bibinfo{volume}{230}, \bibinfo{pages}{6420--6437}.
\bibitem[{Liu et~al.(2011)Liu, Wang and Fox}]{Liu2011}
\bibinfo{author}{Liu, H.}, \bibinfo{author}{Wang, Z.}, \bibinfo{author}{Fox,
  R.O.}, \bibinfo{year}{2011}.
\newblock \bibinfo{title}{A level set approach for dilute non-collisional
  fluid-particle flows}.
\newblock \bibinfo{journal}{Journal of Computational Physics}
  \bibinfo{volume}{230}, \bibinfo{pages}{920--936}.
\bibitem[{Mellet and Vasseur(2007)}]{Mellet2007}
\bibinfo{author}{Mellet, A.}, \bibinfo{author}{Vasseur, A.},
  \bibinfo{year}{2007}.
\newblock \bibinfo{title}{Global weak solutions for a
  vlasov--fokker--planck/navier--stokes system of equations}.
\newblock \bibinfo{journal}{Mathematical Models and Methods in Applied
  Sciences} \bibinfo{volume}{17}, \bibinfo{pages}{1039--1063}.
\bibitem[{Mellet and Vasseur(2008)}]{Mellet2008}
\bibinfo{author}{Mellet, A.}, \bibinfo{author}{Vasseur, A.},
  \bibinfo{year}{2008}.
\newblock \bibinfo{title}{Asymptotic analysis for a
  vlasov-fokker-planck/compressible navier-stokes system of equations}.
\newblock \bibinfo{journal}{Communications in Mathematical Physics}
  \bibinfo{volume}{281}, \bibinfo{pages}{573--596}.
\bibitem[{Patankar and Joseph(2001a)}]{Patankar2001IJMFb}
\bibinfo{author}{Patankar, N.}, \bibinfo{author}{Joseph, D.},
  \bibinfo{year}{2001}a.
\newblock \bibinfo{title}{Lagrangian numerical simulation of particulate
  flows}.
\newblock \bibinfo{journal}{International Journal of Multiphase Flow}
  \bibinfo{volume}{27}, \bibinfo{pages}{1685--1706}.
\bibitem[{Patankar and Joseph(2001b)}]{Patankar2001IJMFa}
\bibinfo{author}{Patankar, N.}, \bibinfo{author}{Joseph, D.},
  \bibinfo{year}{2001}b.
\newblock \bibinfo{title}{Modeling and numerical simulation of particulate
  flows by the eulerian--lagrangian approach}.
\newblock \bibinfo{journal}{International Journal of Multiphase Flow}
  \bibinfo{volume}{27}, \bibinfo{pages}{1659--1684}.
\bibitem[{Potapczuk and Tsao(2019)}]{Potapczuk2019}
\bibinfo{author}{Potapczuk, M.}, \bibinfo{author}{Tsao, J.},
  \bibinfo{year}{2019}.
\newblock \bibinfo{title}{{The Influence of SLD Drop Size Distributions on Ice
  Accretion in the NASA Icing Research Tunnel}}.
\newblock \bibinfo{publisher}{SAE Technical Paper}.
\bibitem[{Prosperetti and Tryggvason(2007)}]{Prosperetti2007}
\bibinfo{author}{Prosperetti, A.}, \bibinfo{author}{Tryggvason, G.},
  \bibinfo{year}{2007}.
\newblock \bibinfo{title}{Computational Methods for Multiphase Flows}.
\newblock \bibinfo{publisher}{Cambridge University Press},
  \bibinfo{address}{Cambridge, UK}.
\bibitem[{Williams(1985)}]{Williams1985}
\bibinfo{author}{Williams, F.A.}, \bibinfo{year}{1985}.
\newblock \bibinfo{title}{Combustion Theory, 2nd ed.}
\newblock \bibinfo{publisher}{Benjamin Cummings}, \bibinfo{address}{Menlo Park,
  CA}.
\bibitem[{Yu(2013)}]{Yu2013}
\bibinfo{author}{Yu, C.}, \bibinfo{year}{2013}.
\newblock \bibinfo{title}{Global weak solutions to the incompressible
  navier–stokes–vlasov equations}.
\newblock \bibinfo{journal}{Journal de Mathématiques Pures et Appliquées}
  \bibinfo{volume}{100}, \bibinfo{pages}{275--293}.

\end{thebibliography}
	
\end{document}